\documentclass[psamsfonts,reqno]{amsart}

\usepackage{amssymb,amsfonts}
\usepackage[all,arc]{xy}
\usepackage{enumerate}
\usepackage{mathrsfs}
\usepackage{xypic}

\usepackage{amsmath,amssymb,amsfonts,wasysym}   
\usepackage{graphics}                           
\usepackage{amsmath}  
\usepackage{mathtools} 
\usepackage{hyperref}                           
\usepackage{comment}
\usepackage[english]{babel}

\newtheorem{thm}{Theorem}[section]

\newtheorem{prop}[thm]{Proposition}

\theoremstyle{definition}
\newtheorem{defn}[thm]{Definition}

\newtheorem{obs}[thm]{Observation}

\newtheorem{exmp}[thm]{Example}

\theoremstyle{remark}

\newcommand*{\im}{\mathop{}\!\mathrm{Im}}
\newcommand*{\id}{\mathop{}\!\mathrm{id}}
\newcommand*{\dimensione}{\mathop{}\!\mathrm{dim}}
\newcommand*{\degree}{\mathop{}\!\mathrm{deg}}

\newcommand*{\spaned}{\mathop{}\!\mathrm{span}}
\newcommand*{\kernel}{\mathop{}\!\mathrm{Ker}}
\newcommand*{\lie}{\mathop{}\!\mathrm{Lie}}

\newcommand*{\Proof}{\mathop{}\!\mathit{Proof}}

\newcommand*{\val}{\mathop{}\!\mathrm{val}}

\newcommand{\bigslant}[2]{{\raisebox{.2em}{$#1$}\left/\raisebox{-.2em}{$#2$}\right.}}

\makeatletter
\let\c@equation\c@thm
\makeatother
\numberwithin{equation}{section}

\title{Prequantization, Geometric quantization, Corrected Geometric Quantization}

\author{Simone Camosso}

\date{}

\begin{document}
{\renewcommand{\thefootnote}{\fnsymbol{footnote}}
\setcounter{footnote}{1}
\footnotetext{\textbf{e-mail}: r.camosso@alice.it, \textbf{address}: I.I.S. ``G.Soleri--A.Bertoni'' Saluzzo (CN), via Traversa del Quartiere n.2, Italy. 
These notes has been written during STRING MATH 2020 (a virtual Zoom event during the Pandemic of Covid-19)}
\setcounter{footnote}{0}
}
\begin{abstract}
A comparison on some facts concerning the geometric quantization of symplectic manifolds is presented here. Criticism and improvements on the sophisticated theory of geometric quantization are presented touching briefly, all the ``salient points of the theory''. The unfamiliar reader can consider this as a ``soft'' introduction to the topic. 
\end{abstract}

\maketitle

\tableofcontents

\section{Introduction}

The intention of these pages is to offer a ``distillation of ideas'' on the geometric quantization. 
The studies that have led to the birth of this line of interest are principally due to these authors: \cite{GS3}, \cite{GS4}, \cite{K}, \cite{So} and \cite{W3}.
On the contrary, the works that inspired this paper are principally two. The first is a recent article of Carosso \cite{C}, where the author describes in a very detailed way the procedure of the geometric quantization with relevant attention and with professional criticism. The second is the celebrated book of Woodhouse \cite{W}, considered as the basic text for ``beginners''. 

The first pages of this article are a review of basic facts on K\"{a}hler manifolds and classical mechanics (especially from the mathematical point of view). After introducing this fundamental formalism with a lot of examples, the author proceeds the examination of the ``quantization procedure''.  The section three is dedicated to the process of prequantization with the description of all prequantum conditions adopted by different authors. In the section four the process of geometric quantization is described introducing the notion of polarization (real and complex) and examining the particular case of compact K\"{a}hler manifolds. In the section five the square--root bundle is introduced in order to pass to a corrected geometric quantization. BKS pairing, used by different authors, is examined in three cases: in the presence of two real polarizations, one real and one complex and, at the end, two complex polarizations. The last section is dedicated to the Bohr-Sommerfeld  subvarieties with recent developments (see \cite{I}).

In a subsection of the section five, the Schr\"{o}dinger equation is derived for the case of the cotangent bundle in a similar way to \cite{C} with the difference we used a result of Albeverio and Mazzucchi (\cite{AM}) instead the heat kernel.

The implicit question regarding the future of the geometric quantization is legitimate and necessary. We will see, during these pages, different criticisms on this rigorous procedure together different problems in the construction of a quantization space that is sufficiently satisfactory. Furthermore, at the present state of the art, an application for the relativity theory would be desirable but not possible. What emerge  in this direction is that the machinery of geometric quantization seems to be not the rigorous solution to the problem of the unification of the relativity and the quantum theory (at least without future modifications). The problem of the geometric quantization of the general relativity has been recently discussed in \cite{MC1} and \cite{MC2}. An interesting variation is given by the deformation quantization, not treated here, but much more suitable. Another possible way in this direction is the method used by Feynman considering a  ``summation of histories'' or probability amplitudes. For this reason an interesting development can be based on the relation between the geometric quantization and the Feynman integral. A seed of this idea can be found in \cite{Si} where the author shows how the Feynman integral can be view as a particular case of BKS pairing. This possibility sure is not the end of the history, as other approaches show this is not the only ``via regia'' (the literature regarding other quantization methods is omitted here due to the large amount to be taken into consideration). A last work that we cite is \cite{Cam2} where the author studied the connection between the geometric quantization (GQ) process and the quantum logic (QL). In this optics the geometric quantization can be a considered as a ``machinery'' that produces Hilbert spaces with interesting properties.

\section{K\"{a}hler manifolds, classical mechanics and symmetries}

\subsection{K\"{a}hler manifolds}

A K\"{a}hler manifold is represented by a quadruple $(M,\omega,J,K)$ where $M$ is a complex manifold of complex dimension $n$, $\omega$ is the symplectic structure locally given by:

\begin{equation}
\label{symplectic}
\omega\,=\,i\sum_{j,l=1}^{n}g_{jl}dz_{j}\wedge d\overline{z}_{l},
\end{equation}
where $g$ is the Riemannian metric on $M$, $J$ is the complex structure associated to $M$ and $K$ is the K\"{a}hler potential associated to $\omega$:

$$ \omega\,=\, i\partial\overline{\partial}K.$$

We recall that in a K\"{a}hler manifold $d\omega=0$ then locally $\omega$ admits a potential $K$. Furthermore the metric $g$ is $J$-invariant and 

$$\omega(v,w)\,=\, g(Jv ,w),$$
\noindent
for every $v,w\in TM$. Furthermore on $M$ it is possible to consider also an hermitian metric, denoted by $h$ and defined by:

$$h\,=\, g-i\omega.$$

The complex structure $J:TM\rightarrow TM$ divides the tangent space in a point $m\in M$, into a direct sum $T_{m}M\,=\, T_{m}M\oplus \overline{T_{m}M}$ where:

$$ T_{m}M\,=\, \{v\in T_{m}M: Jv= iv\}, \ \ \overline{T_{m}M}\,=\, \{v\in T_{m}M: Jv=-iv\}.$$

The first $T_{m}M=\spaned\langle \frac{\partial}{\partial_{z_{j}}}\rangle_{j=1,\ldots,n}$ and the second $\overline{T_{m}M}=\spaned\langle \frac{\partial}{\partial_{\overline{z}_{j}}}\rangle_{j=1,\ldots,n}$. Sometimes these operators are simply denoted by $\partial= \left(\partial_{z_{1}}, \ldots, \partial_{z_{n}} \right)$ and $\overline{\partial}=\left(\partial_{\overline{z}_{1}}, \ldots, \partial_{\overline{z}_{n}} \right)$. Passing through the $2n$ real coordinates denoted by $\{x_{j},y_{j}\}$ (for $j=1, \ldots, n$), we have:

$$\frac{\partial}{\partial_{z_{j}}}\,=\, \frac{1}{2}\left(\frac{\partial}{\partial x_{j}}-i\frac{\partial}{\partial y_{j}}\right), \ \ \frac{\partial}{\partial \overline{z}_{j}}\,=\, \frac{1}{2}\left(\frac{\partial}{\partial x_{j}}+i\frac{\partial}{\partial y_{j}}\right)$$
\noindent
and respectively, the dual basis is $dz_{j}\,=\, dx_{j}+idy_{j}$ and $d\overline{z}_{j}\,=\, dx-idy_{j}$. 

\begin{exmp}
An example of K\"{a}hler manifold is $\left(\mathbb{C}^{n},\omega,J,K\right)$, where:

$$ \omega\,=\, \frac{i}{2}\sum_{j=1}^{n}dz_{j}\wedge d\overline{z}_{j},$$
\noindent 
and the K\"{a}hler potential is given by$K(z)\,=\, \sum_{j=1}^{n} z_{j}\overline{z}_{j}\,=\, \|z\|^{2}$. Here the symplectic form is the standard form $\omega_{0}\,=\, \sum_{j=1}^{n}dx_{j}\wedge dy_{j}$.
\end{exmp}	

\begin{exmp}
	Another example of K\"{a}hler manifold is $\left(\mathbb{PC}^{n}=\bigslant{\mathbb{C}^{n}\setminus \{0\}}{\mathbb{C}^*},\omega_{\text{FS}},J,K\right)$, where:
	
	$$ \omega_{\text{FS}}\,=\, \frac{i}{2}\partial\overline{\partial}\log{\left(1+\|z\|^2\right)},$$
	\noindent 
	is the Fubini--Study form and the K\"{a}hler potential is $K(z)\,=\, \log{\left(1+\|z\|^2\right)}$ (here $\log$ denotes the natural logarithm). 
\end{exmp}	

\begin{exmp}
	The unit disk $\mathbb{D}=\{z\in\mathbb{C}: |z|^2<1\}$ has the structure of K\"{a}hler manifold $\left(\mathbb{D},\omega_{\text{HY}},J,K\right)$ where:
	
	$$ \omega_{\text{HY}}\,=\, \frac{i}{\pi}\frac{dz\wedge d\overline{z}}{\left(1-|z|^2\right)^2}$$
	\noindent 
	is the hyperbolic form and $K(z)\,=\, \frac{1}{\pi}\log{\left(1-|z|^2\right)}$.
\end{exmp}

\begin{exmp}
	A last example is the complex torus $\mathbb{T}=\bigslant{\mathbb{C}}{\Lambda}$ where $\Lambda$ is a lattice on $\mathbb{C}$.  $\left(\mathbb{T},\omega_{\text{flat}},J,K\right)$ has the structure of K\"{a}hler manifold  where:
	
	$$ \omega_{\text{flat}}\,=\, \frac{i}{2}\partial\overline{\partial}H$$
	\noindent 
	and $H:\mathbb{C}\times\mathbb{C}\rightarrow \mathbb{C}$ is a complex map linear in the first factor and respectively antilinear in the second factor such that $H(z,z)\geq 0$ and $H(z,z)=0$ if and only if $z=0$. In this case $K=H$. 
\end{exmp}

\subsection{K\"{a}hler manifolds  and classical mechanics}

In this section we recall the relation between K\"{a}hler manifolds and classical mechanics. Let $(M,\omega,J,K)$ be a K\"{a}hler manifold with $\dim_{\mathbb{C}}{M}=n$ and $\omega$ its symplectic form, then it is possible to interpret $M$ as a model of a classical system. For example $M=T^{\vee}C$ the cotangent space of a configuration space $C$ of dimension $n$. So in this case $\omega$ is given in local symplectic coordinates and classical observables are smooth functions $f$ on $M$.  
To each observable $f\in C^{\infty}(M)$ there is an Hamiltonian vector field $X_{f}\in TM$ defined by the equation:

\begin{equation}
\label{classical_m}
df\,=\, X_{f}\, \lrcorner\,\omega,
\end{equation}
\noindent 
we underline that this notation is equivalent to other common notations $df\,=\, \omega(X_{f},\cdot )$ or $df\,=\, i_{X_{f}}\omega$. In local coordinates:

\begin{equation}
\label{classical_m2}
X_{f}\,=\,\partial_{p_{j}}\partial_{q_{j}} -\partial_{q_{j}}\partial_{p_{j}}.
\end{equation}

For any two observables $f,g$ we can define the Poisson bracket $\{f,g\}$ as:

\begin{equation}
\label{poisson}
\{f,g\}\,=\,X_{g}(f)\,=L_{X_{g}}(f)\,=\,df(X_{g})\,=\, \omega\left(X_{f},X_{g}\right)\,=\, -\{g,f\},
\end{equation}
\noindent 
where $L_{X_{g}}$ is the Lie derivative along $X_{g}\in TM$. We recall that the Lie derivative is defined by the flow $\phi_{t}^{X_{g}}(m): \mathbb{R}\times M\rightarrow M$:

\begin{equation}
\label{lie}
L_{X_{g}}(f)\,=\, \frac{d}{dt} f(\phi_{t}^{X_{g}}(m)),
\end{equation}
\noindent 
where in the previous formula generally there is a tensorial field instead of $f$.

An observable is conserved when $L_{X_{g}}(f)\,=\,0$. In particular there is a smooth function $\mathcal{H}:M\rightarrow \mathbb{R}$, called Hamiltonian, such that determines the trajectories of a classical system via the Hamilton's equations:

\begin{equation}
\label{ODE_flow}
\frac{d}{dt}\phi_{t}(m)\,=\, X_{\mathcal{H}}\phi_{t}(m),
\end{equation}   
\noindent 
for a point $m\in M$. 

The system $(M,\omega,X_{\mathcal{H}})$ is called Hamiltonian system. If $\phi_{t}^{X_{\mathcal{H}}}$ is an integral curve of $X_{\mathcal{H}}$ then the energy function $\mathcal{H}\left(\phi_{t}^{X_{\mathcal{H}}}\right)$ is constant for all $t$ and $\mathcal{H}\left(\phi_{t}^{X_{\mathcal{H}}}\right)\,=\, \mathcal{H}$.
		
Note that $L_{X_{\mathcal{H}}}(\omega)= X_{\mathcal{H}}\,\lrcorner\,d\omega+d(X_{\mathcal{H}}\,\lrcorner\,\omega) = 0+dd\mathcal{H}=0$.
		
Further information on Hamiltonian mechanics can be found in \cite{H}.

From the fact that $\omega$ is closed we have that locally admits a potential $\theta\,=\,\sum_{j=1}^{n}p_{j}dq_{j}$. From the potential we can define the lagrangian associate to an observable $f$ defined as:

\begin{equation}
\label{lagrangian}
\mathcal{L}\,=\, X_{f}\,\lrcorner\,\theta - f.
\end{equation}

The lagrangian mechanics is associated to lagrangian submanifolds of $M$. When the function $f$ is the Hamiltonian $\mathcal{H}$ then $\mathcal{L}\,=\,X_{\mathcal{H}}\,\lrcorner\,\theta - \mathcal{H}$ is the Legendre transform formula and by the Cartan's formula we have that:

$$ L_{X_{\mathcal{H}}}\theta\,=\, X_{\mathcal{H}}\,\lrcorner\,d\theta + d\left(X_{\mathcal{H}}\,\lrcorner\,\theta\right)\,=\, d\mathcal{L}.$$

We recall that $X_{\mathcal{H}}$ generates a diffeomorphism $\rho:M\rightarrow M$ that preserve $\omega$. From the fact that $ \omega\,=\, \rho^{*}\omega$ we have that $d\theta -\rho^{*}d\theta \,=\, 0$ so $\theta -\rho^{*}\theta$ is closed and locally exact. Thus $\theta -\rho^{*}\theta\,=\, dS$, where $S$ is a function on $M=T^{\vee}C$. Let $\gamma$ be an integral curve of $X_{\mathcal{H}}$ with parametrization $t$, then we have that:

$$  (\theta(X_{\mathcal{H}})\circ\gamma)(w)|^{t}_{0}\,=\, (L\circ\gamma)(w)|^{t}_{0} $$
\noindent 
and

$$(\theta(X_{\mathcal{H}})\circ\gamma)(w)|^{t}_{0}\,=\, (\rho_{w}^{*}\theta)(\rho^{-1}_{w^{*}}X_{\mathcal{H}})(\gamma(0))|^{t}_{0}\,=\, (\rho_{w}^{*}\theta)(X_{\mathcal{H}})(\gamma(0))|^{t}_{0}$$
$$\,=\, \theta(X_{\mathcal{H}})(\gamma(0))|^{t}_{0}-(dS(X_{\mathcal{H}})\circ\gamma)(w)|^{t}_{0}$$
\noindent 
because the pushforward of $X_{\mathcal{H}}$ along itself is the identity and, the form $\theta-\rho_{w}^{*}\theta$ is exact. By integration:

\begin{equation}
\label{integral}
S\circ\rho_{w} \,=\, -\int_{0}^{t}\mathcal{L}\circ \gamma (w) dw + c,
\end{equation}
\noindent
where $S\circ \rho_{w}$ is a local phase function (because $X_{\mathcal{H}}$ is complete). We obtain thus a function $\mathcal{S}$ called the generating function:

\begin{equation}
\label{integral}
\mathcal{S}(q,q',t) \,=\, -\int_{0}^{t}\mathcal{L}\circ \gamma (w) dw + c,
\end{equation}
\noindent 
where $\gamma(t)\,=\, \rho_{t}(q,p)\,=\, (q',p')$ and $q(0)=q$ and $q'=q(t)$. The function $\mathcal{S}$ generates a submanifold $\Lambda \subset M$ where $\omega=0$ and where $\dim(M)=\dfrac{n}{2}$. Such manifolds $\Lambda$ are called lagrangian submanifolds. The fiber coordinates $p_{j}$ on $\Lambda$ are given by:
\begin{equation}
\label{hamilton_J}
p_{j}\,=\, \frac{\partial \mathcal{S}}{\partial q_{j}},
\end{equation}
\noindent 
for $j=1,\ldots , n$. In order to pass from $S(q,p(q'),t)$ to $\mathcal{S}(q,q',t)$ we use the equations $(\ref{hamilton_J})$.

Details on lagrangian and hamiltonian mechanics are in \cite{C} and \cite{W}.

\subsection{K\"{a}hler manifolds  and classical mechanics with symmetries}

Let $G$ be a finite dimensional compact Lie group and $\mu^{G}:G\times M\rightarrow M$ an action on the K\"{a}hler manifold $(M,\omega,J,K)$. Let $\mathfrak{g}=\lie{(G)}$ be the Lie algebra associated to $G$ and $\mathfrak{g}^{\vee}$ the dual space. If we have an $n$--form $\omega$ on $M$ it is associated a map $\Psi:M \rightarrow \bigwedge^{n}\mathfrak{g}^{\vee}$ defined as:

\begin{equation}
\label{moment1}
\Psi(m)\,=\, \mu_{m}^{G\vee}\omega,
\end{equation}
\noindent 
where $\mu_{m}^{G}:G\rightarrow M$ is $\mu_{m}^{G}(g)\,=\, g\cdot m$ (observe that if the action of $G$ is symplectic then $\mu_{m}^{G\vee}\omega\,=\, \omega$). Furthermore if $\omega$ is closed, then also $\Psi$ ($\forall m\in M$) and $\Psi:M\rightarrow Z^{2}(\mathfrak{g})$. If in addition $H^{1}(\mathfrak{g})\,=\, H^{2}(\mathfrak{g}) \,=\, \{0\}$ then $\forall m\in M$ there is a unique element $\Phi(m)$ in $\mathfrak{g}^{\vee}$ such that:

\begin{equation}
\label{moment1}
d\Phi(m)\,=\, \Psi(m).
\end{equation}
There is a unique map $\Phi:M\rightarrow \mathfrak{g}$ such that:

\begin{equation}
\label{moment2}
d\Phi\,=\, \Psi.
\end{equation}

We can be more explicit with $\Phi$. In fact, under our assumptions, there is an homeomorphism between the Lie algebra $\mathfrak{g}$ and the hamiltonian vector fields on $M$ (see \cite{GS1}). In particular there is a unique homeomorphism $\lambda:\mathfrak{g}\rightarrow C^{\infty}(M)$ such that $\lambda(\xi)=f_{\xi}$ where:

\begin{equation}
\label{moment3}
\xi_{M} \,\lrcorner\,\omega \,=\, df_{\xi}.
\end{equation}

The map $f_{\xi}$ depends linearly by $\xi$ so $\forall m\in M$ we can consider $\Phi(m)\in\mathfrak{g}^{\vee}$ given by:

\begin{equation}
\label{moment4}
\langle\Phi(m),\xi \rangle \,=\,  f_{\xi}(m),
\end{equation}
\noindent 
where $\langle\cdot, \cdot \rangle$ is the pairing between $\mathfrak{g}$ and $\mathfrak{g}^{\vee}$. In other words:

\begin{equation}
\label{moment5}
d\langle \Phi(m),\xi \rangle \,=\,\xi_{M} \,\lrcorner\,\omega,
\end{equation}
\noindent 
for every $\xi \in\mathfrak{g}$ where $\Phi: M \rightarrow \mathfrak{g}^{\vee}$ is called the moment map. The previous equation can be written in the following form:

\begin{equation}
\label{moment6}
\langle d_{m}\Phi(v),\xi \rangle \,=\,\omega_{m}\left(\xi_{M}(m),v\right),
\end{equation}
\noindent 
where $d_{m}\Phi:T_{m}M \rightarrow \mathfrak{g}^{\vee}$ for every $m\in M$ and $v\in T_{m}M$.

The moment map has different properties. The first is that $d_{m}\Phi$ is the transpose map of the valuation $\val : \mathfrak{g}\rightarrow T_{m}M$. We have that:

\begin{equation}
\label{moment7}\ker{(d_{m}\Phi)}\,=\, \mathfrak{g}_{M}(m)^{\perp},
\end{equation} 
\noindent
where $\mathfrak{g}_{M}(m)$ is the subspace of vectors $\xi_{M}(m)$. We have that $G$ acts on $M$ transitively so $T_{m}M\,=\,\mathfrak{g}_{M}(m)$ for every $m\in M$. We have also that:

\begin{equation}
\label{moment8} \im{(d_{m}\Phi)^{0}}\,=\, \{\xi\in\mathfrak{g}: \xi_{M}(m)=0\},
\end{equation} 
\noindent 
is the stabilizator of $m\in M$. The stabilizator is discrete if and only if $d_{m}\Phi$ is surjective. The last property is a conseguence of the fact that $G$ acts on $\mathfrak{g}$ by the adjoint representation, $\val$ is a $G$--morphism, thus $\Phi$ is a $G$--morphism:

\begin{equation}
\label{moment9} \Phi(g\cdot m)=g\cdot \Phi(m),
\end{equation} 
\noindent 
for every $g\in G$ and $m\in M$. 

\begin{exmp}
	Let us consider $\mathbb{R}^{2}$ with $\omega=d\,x\wedge d\,y=rd\,r\wedge d\,\theta$. Let us consider as $G$ the circle $S^{1}$ that acts with rotations. The generator of the action is the field $x\frac{\partial}{\partial y}-y\frac{\partial}{\partial x}=\frac{\partial}{\partial \theta}$. Thus $\Phi^{S^{1}}=-\frac{1}{2}(x^{2}+y^{2})=\frac{1}{2}r^{2}$ is the moment map. 
	
	In general on $\left(\mathbb{C}^{n}\equiv \mathbb{R}^{2n},\omega\right)$ where $\omega=\sum_{i=1}^{n}d\,x_{i}\wedge d\,y_{i}$, with the action of $S^{1}$ and heights $m_{1}, \ldots, m_{n}\in \mathbb{Z}$, the action has generator the vector field $\sum_{i=1}^{n} m_{i}\frac{\partial}{\partial \theta_{i}}$ and moment map given by:
	
	$$ \Phi^{S^{1}}(x_{1}, \ldots,x_{n},y_{1},\ldots,y_{n})=-\frac{1}{2}\sum_{i=1}^{n}m_{i}\left(x_{i}^{2}+y_{i}^{2}\right).$$
\end{exmp}

\section{Prequantization}

\subsection{Prequatization conditions}

Recalling that a closed surface $\Sigma$ of $M$ is a surface that is compact and without boundary, we state the prequantum condition \textbf{PC1} as the following.

\begin{center}\textbf{PC1} The integral of $\omega$ over any closed 2--surface $\Sigma \subset M$ is an integral multiple of $2\pi \hbar$.\end{center}

This prequantum condition \textbf{PC1} is necessary for the existence of the hermitian line bundle $L$ over $M$. In the case of $M$ simply connected \textbf{PC1} is also sufficient. Assume $M$ simply connected and $m_{0}\in M$ be a base point. Let us consider the following set:

$$\{(m,z,\gamma): m\in M, z\in\mathbb{C}, \gamma \ \ \text{a piecewise smooth path from} \ \ m_{0} \ \ \text{to} \ \ m\}$$
\noindent
and the equivalence relation $\sim$ defined as:

$$(m,z,\gamma)\sim (m',z',\gamma') \Leftrightarrow m=m' \wedge z'\,=\, ze^{\frac{i}{\hbar}\int_{\Sigma_{1}}\omega},$$
\noindent 
where $\Sigma_{1}$ is any oriented 2--surface with boundary made up of $\gamma$ (from $m_{0}$ to $m$) and $\gamma'$ (from $m$ to $m_{0}$). The surface $\Sigma_{1}$ exists because $M$ is simply connected. We define the line bundle $L$ as: 

$$L\,=\, \bigslant{\{(m,z,\gamma): m\in M, z\in\mathbb{C}, \gamma \ \ \text{a piecewise smooth path from} \ \ m_{0} \ \ \text{to} \ \ m\}}{\sim}.$$

We can define operations of addition and scalar multiplication between the fibres:

$$ [(m,z,\gamma)]+[(m,z',\gamma)]\,=\, [(m,z+z',\gamma)], \ \ k[(m,z,\gamma)]\,=\, [(m,kz,\gamma)],$$
\noindent 
with $k\in\mathbb{C}$ and $[(m,z,\gamma)],[(m,z',\gamma)]\in L$. Trivializations of $L$ are determined locally by symplectic potentials. In fact let us assume that $\theta_{j}$  is a symplectic potential on a simply connected open set $U_{j}\subset M$ of a collection $\{U_{j}\}_{j\in J}$. Let us consider a point $m_{1}\in U_{j}$ and a curve $\gamma_{0}$ from $m_{0}$ to $m_{1}$. We define locally a section $s$ of $L$ in $U_{j}$ by

$$s_{j}(m)\,=\, \left[\left(m,e^{-\frac{i}{\hbar}\int_{\gamma_{1}}\theta_{j}},\gamma\right)\right],$$
\noindent
where $\gamma_{1}$ is any curve from $m_{1}$ to $m$ in $U_{j}$ and $\gamma$ is the curve from $m_{0}$ to $m$ obtained from $\gamma_{0}$ and $\gamma_{1}$. We observe that a different choice of $\gamma_{1}$ gives the same value of $s_{j}(m)$ and that a different choice of $m_{1}$ or $\gamma_{0}$ gives the same section multiplied by a constant of modulus one. The effect of replacing $\theta_{j}$ by $\theta_{k}\,=\, \theta_{j}+du_{kj}$, with $u_{kj}(m_{1})=0$, is $s_{k}(m)=e^{-\frac{i}{\hbar}u_{kj}}s_{j}(m)$. 

Now we can assume that the line bundle $L\rightarrow M$ exists. Let us consider the parallel transport of a section $s$ respect to $\nabla$ around a loop $\gamma$ of $m\in M$. Assume that $\gamma$ is the boundary of a $2$--surface $\Sigma_{1}$ contained in the domain of $\theta$ a symplectic potential. Solving the parallel transport equation $\nabla_{\dot{\gamma}}s=0$, the result is equivalent to a linear transformation $L_{m}\rightarrow L_{m}$ given by the multiplication of $e^{\frac{i}{\hbar}\int_{\gamma}\theta}$ that, by the Stokes' theorem, is equivalent to $e^{\frac{i}{\hbar}\int_{\Sigma_{1}}\omega}$. Now let us consider a second surface $\Sigma_{2}$ with boundary $\gamma$ such that $\Sigma\,=\, \Sigma_{1}\cup \Sigma_{2}$ is a closed $2$-surface in $M$ (in other terms gluing together $\Sigma_{1}$ and $\Sigma_{2}$ in $\gamma$ we obtain the closed 2--surface $\Sigma$). In a similar way the parallel transport gives a linear transformation by the multiplication of $e^{-\frac{i}{\hbar}\int_{\Sigma_{2}}\omega}$, where the minus is because the boundary of $\Sigma_{2}$ is $-\gamma$. By the uniqueness of the solution of the differential equation associated to the parallel trasport we have that:

$$e^{\frac{i}{\hbar}\int_{\Sigma_{1}}\omega}\,=\, e^{-\frac{i}{\hbar}\int_{\Sigma_{2}}\omega},$$
\noindent
that is

$$ e^{\frac{i}{\hbar}\int_{\Sigma}\omega}\,=\, 1,$$
\noindent
the last equation is equivalent to \textbf{PC1} because $\int_{\Sigma}\omega\in 2\pi \hbar\mathbb{Z}$.

The prequantum condition is related to the existence of the hermitian line bundle $L$ over $M$ also through the Weil theorem.

\begin{thm}[Weil, 1958, \cite{We}]
Let $M$ be a smooth manifold and $\omega$ a real, closed 2--form whose cohomology class $\left[\frac{\omega}{2\pi\hbar}\right]$ is integral. Then there is a unique hermitian line bundle $L$ over $M$ with unitary connection $\nabla$ so that $c_{1}(L)=\left[\frac{\omega}{2\pi\hbar}\right]$. 
\end{thm}

The converse is also true. In fact let $\{U_{j}\}$ be a contractible open cover of $M$. By assumption there is a collection of $1$--forms $\theta_{j}\in \Omega^{1}(U_{j})$ such that $\frac{\omega}{2\pi \hbar}\,=\, d\theta_{j}$ on $U_{j}$. We can find $u_{jk}\in C^{\infty}(U_{j}\cap U_{k})$ such that $du_{jk}\,=\, \theta_{j}-\theta_{k}$ whenever $U_{j}\cap U_{k}\,\not=\, \emptyset$.
On $U_{j}\cap U_{k}\cap U_{l}$ (whenever $U_{j}\cap U_{k}\cap U_{l}\,\not=\,\emptyset$) we have that $ du_{jk}+du_{kl}+du_{lj}\,=\, 0$, so the function $u_{jkl}\,=\,u_{jk}+u_{kl}+u_{lj}$ is constant. They define a C\v{e}ch cohomology class in $H^{2}(M,\mathbb{R})$. 
Now $\theta_{j}-\theta_{k}\,=\, d\log{g_{jk}}$ whenever $U_{j}\cap U_{k}\,\not=\, \emptyset$ and, the de Rham isomorphism send $c_{1}(L)$ to $\left[\frac{\omega}{2\pi \hbar}\right]\in H^{2}(M,\mathbb{R})$ with:

$$ u_{jkl}\,=\, \frac{1}{2\pi i\hbar}\left(\log{g_{jk}}+\log{g_{kl}}+\log{g_{lj}}\right).$$
 
By assumption $g_{kl}$ are transition functions and must satisfy the cocycle condition: 

$$ g_{jk}g_{kl}g_{lj}\,=\, 1,$$
\noindent 
that is

$$ e^{2\pi i\hbar u_{jkl}}\,=\,1.$$

So $\left[\frac{\omega}{2\pi \hbar}\right]\in H^2(M,\mathbb{Z})$. 

We have a second version of the prequantization condition \textbf{PC2}.

\begin{center}\textbf{PC2} $\left[\frac{\omega}{2\pi \hbar}\right]\in H^2(M,\mathbb{Z})$. \end{center}

Note that in \textbf{PC2} it is not required $M$ to be simply connected, so it is more general than \textbf{PC1}.

Another prequantization of a symplectic manifold $(M,\omega)$ consist in a $S^1$--bundle $(X,\alpha)$ with the projection $\pi:X \rightarrow M$ and $\alpha$ an $S^{1}$ invariant 1--form such that satisfy the prequantization condition \textbf{PC3}.

\begin{center}\textbf{PC3} $d\alpha\,=\, \frac{1}{2\pi \hbar}\pi^*\omega$. \end{center}

The $S^1$--bundle $(X,\alpha)$ it is also called the circle bundle and is defined as:

$$ X\,=\, \{(m,\lambda):m\in M,\lambda\in L^{\vee}_{m}, h(\lambda,\lambda)=1\}\subset L^{\vee},$$
\noindent 
where $h$ is the hermitian metric on $L$. We have that $(X,\alpha)$ is also called a contact manifold. 

The last way to state the prequantization condition derive from \textbf{PC3} and consists relating the curvature form of the connection $\nabla$ of the line bundle $L$ and the symplectic form $\omega$. In other terms if $\Theta$ is the curvature form of the unique covariant derivative $\nabla$ on $L$ compatible with both the complex and hermitian structures, we have the prequantization condition \textbf{PC4}.

\begin{center}\textbf{PC4} $\Theta\,=\,-\frac{i}{\hbar}\omega$. \end{center}

\subsection{Examples}

\begin{exmp}
Let $H:\mathbb{C}\times \mathbb{C} \rightarrow \mathbb{C}$ be a complex linear map in the first factor and antilinear in the second. Assume also that $H(z,z)\geq 0$ and that is zero if and only if $v=0$. Let us define:

\begin{equation}
\label{toriquant11}
\omega\,=\, \frac{i}{2}\partial\overline{\partial} H,
\end{equation}
\noindent 
then $\omega$ is a K\"{a}hler form for the complex torus $\mathbb{T}\,=\, \bigslant{\mathbb{C}}{\Lambda}$ that is invariant under the action of the lattice $\Lambda$. Furthermore it is possible to prove that the torus $\mathbb{T}$ is quantizable if and only if $\im{H}(\Lambda\times \Lambda) \subset \mathbb{Z}$. This condition it is in fact equivalent to the integrability condition \textbf{PC4}. We can try to see this in one direction. Let us assume  $\im{H}(\Lambda\times \Lambda) \subset \mathbb{Z}$, then the image of $\Lambda\times \Lambda$ through $H$ is in $\mathbb{Z}$. This condition ensure the existence of a complex line bundle $L=L(\chi,H)\cong \mathbb{C}\times \mathbb{T}$ where $\chi$ is a semicharacter associated to $H$. On $L$ we can define an hermitian structure defined as:

\begin{equation}
\label{hermitiantoruslinebundle}
h(\theta(z),\theta(z))\,=\, F(z)|\theta(z)|^2,
\end{equation}
\noindent 
where $\theta(v)$ is the holomorphic section of $L$, also called theta function such that:

\begin{equation}
\label{thetaeq1}
\theta(z+\lambda)\,=\,A(\lambda,z)\theta(z),
\end{equation} 
\noindent 
and $A$ is a factor of automorphy on $\Lambda\times \mathbb{C}$ (a good reference is \cite{LB}). The second term in the definition of $h$ is $F(z)$ that is a map $F:\mathbb{C}\rightarrow \mathbb{R}_{+}$ defined as:

\begin{equation}
\label{thetaeq2}
F(z)\,=\,e^{-\pi H(z,z)}.
\end{equation} 

It is possible to see that with this definition, the function $h$ is invariant under the action of $\Lambda$ and $h$ defines an hermitian structure on $L$. Now analyzing the curvature form of the line bundle $L$ we find that it is equal to $-\partial\overline{\partial}\log{h}=-\partial\overline{\partial}\log{F}=\pi\partial\overline{\partial}H$. This last term is equal to $-2\pi i \omega$ that is the condition \textbf{PC4} with the Planck constant equal to $1$. 
The converse of the proof is proved in \cite{LB}.
\end{exmp}

\begin{exmp}
Let $M=S^2$ be the sphere in $\mathbb{R}^{3}$ of radius $r$. The sphere is a K\"{a}hler manifold with symplectic structure:

\begin{equation}
\label{spherequant11}
\omega\,=\, r\sin{\varphi}d\varphi\wedge d\vartheta,
\end{equation}
\noindent 
that is a 2--form on $S^2$ where, $\vartheta\in[0,2\pi)$, $\varphi\in [0,\pi]$ and $r\in [0,+\infty)$. In order to be quantizable the integral:

$$ \int_{S^{2}}\omega \,=\, 4\pi r,$$
\noindent 
must be equal to $2\pi \hbar n$ for some $n\in\mathbb{Z}$. We conclude that not all spheres satisfy the prequantization condition \textbf{PC1} but only the spheres with $r=\frac{n}{2}\hbar$.
\end{exmp}

\begin{exmp}
	Let $M= S^2$ be the sphere realized as section of the light-cone by the hyperplane $x^0=ct=1$. It is show in \cite{V} that $M$ can be quantized in the  `` sense of Souriau''. The quantization is obtained as $S^1$--fiber bundle $(S^3,\pi)$ where $\pi$  are the so called ``KS--transformations'' (Kustaankeimo and Stiefel transformations for the regularization of the Kepler problem) which associates to each vector in the light--cone a one--index spinor. In practice, in \cite{V}, it is described how from the KS--trasformations we obtain the Hopf fibering of the sphere $S^3$. In this case the quantization condition is the same of the previous example for $r=1$:
	
	$$ \int_{S^{2}}\omega \,=\, 4\pi .$$
\end{exmp}

\begin{exmp}
Let $M=\mathbb{D}$ be the unit disk. It is a K\"{a}hler manifold with the hyperbolic form and, the prequantization condition $\textbf{PC1}$, is trivially satisfied. The line bundle is trivial $L=\mathbb{D}\times \mathbb{C}\rightarrow \mathbb{D}$  and $h(z,t)\,=\, \left(1-|z|^{2}\right)^2|t|^2$ is the hermitian structure with $(z,t)\in \mathbb{D}\times \mathbb{C}$. An analogue argument shows that the complex space $\mathbb{C}^{n}$ is prequantizable, in the same way.
\end{exmp}

\begin{exmp}
	Let $M=\Sigma$ be a compact Riemann surface of genus $g\geq 2$. We can think the Riemann surface as the quotient $\bigslant{\mathbb{D}}{G}$, where $\mathbb{D}$ is the unit disk of $\mathbb{C}$ and $G$ is the subgroup of $SU(1,1)$ of fractional linear transformation. An element $g\in G$ is represented by a matrix:
	
	\begin{equation}
	\label{riemannsurff}
	g\,=\, \left(\begin{array}{cc} a &  b \\ \overline{b} & \overline{a} \end{array}\right),
	\end{equation}
\noindent 
such that $|a|^2-|b|^2\,=\,1$. The action on $\mathbb{D}$ is defined by $g\cdot z \,=\, \frac{az+b}{\overline{b}z+\overline{a}} $. 	
We have that $\Sigma$ with $\omega_{\text{HY}}$ is a K\"{a}hler manifold because $\omega_{\text{HY}}$ is invariant by the action of $G$.

Let us consider $\{(U_{j},z_{j})\}_{j\in J}$ be a complex atlas on $\Sigma$ and $s_{j}(z)$ be holomorphic functions on $U_{j}$. We can define global sections $s(z)\,=\, s_{j}(z)dz_{j}$ such that on $U_{j}\cap U_{k}$ we have that $s_{j}(z)dz_{j} \,=\, s_{k}(z)dz_{k}$. Now from complex analysis we have that:

$$ \frac{dz_{k}}{dz_{j}}\cdot \frac{dz_{l}}{dz_{k}} \,=\, \frac{dz_{l}}{dz_{j}} $$
\noindent 
that is, defining $\gamma_{jk}=\frac{dz_{k}}{dz_{j}}$, equivalent to $\gamma_{jk}\gamma_{kl}\gamma_{lj}=1$ for every $z\in U_{j}\cap U_{k}\cap U_{l}$. Then there exists a canonical bundle $K$. If we consider the projection $p:\mathbb{D} \rightarrow  \bigslant{\mathbb{D}}{G}$ the pull-back $p^*(K)$ is holomorphically trivial and its holomorphic global sections are of the form $s(z)dz$ on $\mathbb{D}$. We can think to global holomorphic sections on $K$ as $1$--forms of type $(1,0)$ on $\mathbb{D}$ invariant under the action of $G$. They are forms $s(z)dz$ such that:

\begin{equation}
\label{condition_invariant_riemann}
s(g(z))\cdot \frac{1}{\left(\overline{b}z+\overline{a}\right)^2}dz\,=\, s(z)dz,
\end{equation}
\noindent 
for every $g\in G$.
We can proceed as in the case of the torus defining an hermitian structure:

\begin{equation}
\label{hermitian}
h(s(z),s(z))\,=\, (1-|z|^2)^2|s(z)|^2,
\end{equation}
\noindent 
that is invariant for every section of $K$ and $g\in G$. 

At the end we see that the curvature form of the line bundle $K$ is equal to $-\partial\overline{\partial}\log{h}=-\partial\overline{\partial}\log{(1-|z|^2)^2}=\frac{2dz \wedge d\overline{z}}{(1-|z|^2)^2}$. This last term is equal to $-2\pi i \omega_{\text{HY}}$ that is the condition \textbf{PC4} with the Planck constant equal to $1$. 
\end{exmp}

\begin{exmp}
	A very remarkable example is the complex projective space $M=\mathbb{PC}^{n}$. In this case we can start to show that $\mathbb{PC}^{1}$ satisfy the condition \textbf{PC2} with the convention that $h=1$ (in this case it is a Riemann surface and the result it is true). To see this we can use these relations:
	
	\begin{equation}
	\label{degree_tauto_chern_class}
	\degree{(L)}\,=\,\int_{\mathbb{PC}^{1}}c_{1}(L)\,=\, \frac{i}{2\pi}\int_{\mathbb{PC}^{1}}\partial\overline{\partial}\log{(1-|z|^2)},
	\end{equation}
	\noindent 
	where in this case $L\,=\, \mathcal{O}(1)$ is the hyperplane bundle. Now we can directly perform the integration:
	
		\begin{equation}
	\label{degree_tauto_chern_class}
	\degree{(\mathcal{O}(1))}\,=\,\frac{1}{\pi}\int_{\mathbb{R}^{2}}\frac{1}{(1+x^2+y^2)^2}dxdy\,=\,  \frac{1}{\pi}\int_{0}^{+\infty}\int_{0}^{2\pi}\frac{r}{(1+r^2)^{2}}drd\vartheta\,=\, 1.
	\end{equation}

Thus we have that $c_{1}(\mathcal{O}(1))\,=\,\left[\omega_{\text{FS}}\right]\in H^{2}(\mathbb{PC}^{1},\mathbb{Z})$, that is \textbf{PC2}. The result is true for general $n$ and to see this it is possible to use the chain of isomorphisms: $H^{2}(\mathbb{PC}^{n},\mathbb{Z})\rightarrow H^{2}(\mathbb{PC}^{n-1},\mathbb{Z}) \rightarrow \ldots \rightarrow H^{2}(\mathbb{PC}^{1},\mathbb{Z})$, in order to have $\left[\omega_{\text{FS}}^n\right]\in H^{2}(\mathbb{PC}^{n},\mathbb{Z})$.

The result that the K\"{a}hler form $\omega_{\text{FS}}$ is integral is also showed in \cite{GH}.
\end{exmp}

\begin{exmp}
	Let $M=S^2\times S^2$ be the product of spheres with the same radius $r$. The manifold $M$ is a K\"{a}hler manifold with symplectic structure:
	
	\begin{equation}
	\label{spherequant1122}
	\omega\,=\, p_{1}^{*}\omega_{S^2}+p_{2}^{*}\omega_{S^2},
	\end{equation}
	\noindent 
	where $p_{1},p_{2}$ are the projections on the two spheres with radius $r\in [0,+\infty)$ and we define the symplectic structure as the pullback respectively of the two K\"{a}hler structures on $S^2$. In order to be quantizable we have that the integral:
	
	$$ \int_{S^{2}\times\{p'\}}\omega \,=\, 4\pi r,$$
	\noindent 
	for every $S^{2}\times\{p'\}\subseteq S^{2}\times S^{2}$, must be equal to $2\pi \hbar n$ for some $n\in\mathbb{Z}$. We conclude that $r=\frac{n}{2}\hbar$, in this case the product of the two sphere is prequantizable. We observe also that if the radii of the two spheres are incommensurable then $S^{2}\times S^2$ is not prequantizable!
\end{exmp}

\begin{exmp}
	In this example we describe the prequantization of the Kepler manifold $X$ defined as:
	
	$$ X\,=\,\{(\xi,x)\in \mathbb{R}^{n+1}\times \mathbb{R}^{n+1}:\xi\cdot \xi =1, x\cdot x =0, x\not=0\},$$
	\noindent 
	where $\cdot$ is the usual Euclidean scalar product. 
	Usually, the manifold $X$ is denoted also with $T^{\vee}(S^{n})\setminus \{0\}$. 
	The Kepler manifold has the following symplectic form:
	
	\begin{equation}
	\label{spherequant1122Kepler}
	\omega_{\text{Kepler}}\,=\, \sum_{j=1}^{n+1}dx_{j}\wedge d\xi_{j}.
	\end{equation}
	
	We can think $X$ as a complex manifold with the identification of $X$ with the complex light cone:
	
	$$ C\,=\,\{z\in\mathbb{C}^{n+1}: z\cdot z=0,z\not=0\}.$$
	
	The identification can be realized using the map: $m:X\rightarrow C$ such that $(\xi,x)\mapsto \|x\|\cdot\xi+ix$. We have that $(X,J,\omega)$ is a K\"{a}hler manifold where $J$ is the pullback of the complex structure of $C^{n+1}$ and 
	
	\begin{equation}
	\label{spherequant1122Kepler}
	\omega_{\text{Kepler}}\,=\, 2i\overline{\partial}\partial \sqrt{\frac{z\cdot \overline{z}}{2}}.
	\end{equation}
	
	The K\"{a}hler form is exact, trivially integral and there is a quantum line bundle $L$. Moreover for $n \geq 3$ we have that $X$ is simply connected and $L$ is holomorphically trivial over $X$ (details are in \cite{R}).
\end{exmp}

\begin{exmp}
	Let $M=S^1\times \mathbb{R}$ be the cylinder. It can be identified with $T^{\vee}S^{1}$. It is a K\"{a}hler manifold with $\omega=dp\wedge d\varphi$. There is a symplectic potential $\theta=pd\varphi$ and $\omega=d\theta$ is globally exact and the prequantum line bundle $L$ is trivial (but not unique!). 
\end{exmp}

\section{Geometric quantization}

\subsection{The Dirac axioms and quantum operators}

In his work \cite{D}, P. Dirac defines the quantum Poisson bracket $[\cdot,\cdot]$ of any two variables $u$ and $v$ as:

\begin{equation}
\label{dirac}
uv-vu=i\hbar [u,v]
\end{equation}
\noindent 
where $\hbar$ is the Plank constant $h$ over $2\pi$. The formula $(\ref{dirac})$ is one of the basic postulates of quantum mechanics. 
We can summarize these postulates as follows. To start we fix a symplectic manifold $(M,\omega)$ of dimension $n$, with $\omega$ the corresponding symplectic structure and an Hilbert space $\mathfrak{H}$. The quantization is a ``way'' to pass from the classical system to the quantum system. In this case the classical system (or phase space) is described by the symplectic manifold $M$ and the Poisson algebra of smooth function on $M$ denoted by $(C^{\infty}(M),\{\cdot,\cdot\})$. 
The quantum system is described by $\mathfrak{H}$. We define ``quantization'' a map $Q$ from the subset of the commutative algebra of observables $C^{\infty}(M)$ to the space of operators in $\mathfrak{H}$. Let $f\in C^{\infty}(M)$ be an observable we have that $Q(f):\mathfrak{H}\rightarrow \mathfrak{H}$ is the corresponding quantum operator.
We can summarize the quantum axioms in this scheme:

\begin{itemize}
	\item[1.] Linearity: $Q(\alpha f+\beta g)=\alpha Q(f)+\beta Q(g)$, for every $\alpha,\beta$ scalars and $f,g$ observables;
	\item[2.] Normality: $Q(1)=\id$, where $\id$ is the identity operator;
	\item[3.] Hermiticity: $Q(f)^{*}=Q(f)$;
	\item[4.] (Dirac) quantum condition: $[Q(f),Q(g)]=-i\hbar Q(\{f,g\})$;
	\item[5.] Irreducibility condition: for a given set of observables $\{f_{j}\}_{j\in J}$, with the property that for every other $g\in C^{\infty}(M)$, 
	such that $\{f_{j},g\}=0$ for all $j\in J$, then $g$ is constant. We can associate a set of quantum operators $\{Q(f_{j})\}_{j\in J}$ such that for every other operator $Q$ that commute with all of them is a multiple of the identity.
\end{itemize}

The last postulate states that in the case we consider a connected Lie group $G$ we say that is a group of symmetries of the physical system if we have the two 
following irreducible representation: one as symplectomorphisms acting on $(M,\omega)$ and another as unitary transformations acting on $\mathfrak{H}$. 
For many details about these postulates look \cite{ELRM}.

\begin{exmp}[Schr\"{o}dinger quantization]
	Let $M=\mathbb{R}^{2n}$ and $(q_{j},p_{j})$ be the canonical coordinates of position and momentum. In this case $Q(q_{j})=q_{j}$ that acts as multiplication and 
	$Q(p_{j})=-i \hbar\frac{\partial}{\partial q_{j}}$. The Hilbert space $\mathfrak{H}=L^{2}(\mathbb{R}^{n},dq_{j})$ and there are the following relations of commutations:
	
	$$ [Q(q_{k}),Q(q_{j})]=[Q(p_{k}),Q(p_{j})]=0, \ \ \ \ [Q(q_{k}),Q(p_{j})]=i\hbar\delta_{kj}\id.$$
\end{exmp}

Let us examine the form of these quantum operators. We can start considering the case where $M=T^{\vee}C$ with $C$ the configuration space. Now $M$ is a symplectic space with $\omega$ as symplectic form and is prequantizable with an hermitian line bundle $(L,h)$. Proceeding in the choice of a potential $\theta$ it is possible to construct $Q(f)$ for each observable $f\in C^{\infty}(M)$ satisfying the quantum axioms. This operator has the following form:

\begin{equation}
\label{qoperator1}
Q(f)\,=\, -i\hbar X_{f} - X_{f}\,\lrcorner\,\theta + f.
\end{equation}

These operators $Q_{f}$ can be ``glued'' in order to form a global operator on the sections of the corresponding line bundle $L$. 

In general if we have a K\"{a}hler manifold $M$ that satisfies the prequantum condition \textbf{PC4}, and if $s$ is a section of the line bundle $L$, for every $f\in C^{\infty}(M)$, we have an Hamiltonian vector field $X_{f}$ and the operator $\nabla_{X_{f}}$ that acts on $L^{2}(M,L,dq^{j})$. So if we define the operator:

\begin{equation}
\label{KonstantS}
Q(f)s=-i\hbar\nabla_{X_{f}}s+f\cdot s,
\end{equation}
\noindent 
now it is an hermitian linear operator and, when $f$ is constant, is only the multiplication by $f$. The formula $(\ref{KonstantS})$ satisfy the Dirac postulate 4.

\subsection{K\"{a}hler polarizations}

Polarizations of symplectic manifolds are introduced in order to have a dependence of sections of $L$ (that are waves functions) by half the coordinates of the configuration space. 

A complex K\"{a}hler polarization of $M$ (of complex dimension $n$) is a smooth complex distribution $D$ (subbundle of $TM$) that is a map that to each point $m\in M$ assigns a linear subspace $D_{m}$  of $T_{m}M$, such that:

\begin{itemize}
	\item[1.] $\dimensione{(D_{m})}=n$;
	\item[2.] $D\cap\overline{D}=\emptyset$;
	\item[3.] $(D\oplus\overline{D})_{m}=T_{m}M_{\mathbb{C}}$;
	\item[4.] $D$ is involutive, $[D,D]\subset D$;
	\item[5.] $D_{m}$ is isotropic, that is $\omega|_{D}=0$;
	\item[6.] For every $\xi\not= 0$ in $D$, $i\omega(\xi,\overline{\xi})>0$.
\end{itemize} 

First note that if $D$ is a polarization, also $\overline{D}$ is a polarization. Second by the Frobenius theorem $D$ is integrale: that is for each point $m\in M$ there is an integrable submanifold $N$ (of $M$) whose tangent space at $m$ is $D_{m}$. These $N$ are called leaves of $D$. 

A section $s$ of the line bundle $L$ over $M$ is called polarized if:

\begin{equation}
\label{polarized1}\nabla_{\xi}s\,=\,0,
\end{equation}
\noindent 
for all $\xi\in D$. 

A polarization is real if $D\,=\, \overline{D}$. 

A reference for K\"{a}hler polarization is \cite{RCS}. Furthermore, if a K\"{a}hler manifold has one K\"{a}hler polarization then it has many, that is changing polarization the Hilbert space change. A good tool in order to consider possible relations between these Hilbert spaces is the pairing we will see next. 

\begin{exmp}[K\"{a}hler polarization]
Let $(M,\omega,J,K)$ be a K\"{a}hler manifold, then a K\"{a}hler polarization $D$ consists in a submanifold spanned by vectors $\xi$ such that $\omega|_{D}=0$. For example in the Boson--Fock space of dimension $1$, that is $\mathbb{C}$ with $\omega\,=\, \frac{i}{2}dz\wedge d\overline{z}$ and $K(z)=\frac{1}{2}z\cdot \overline{z}$. The K\"{a}hler polarization $\overline{D}\,=\, \langle\partial_{\overline{z}}\rangle$ is spanned by the antiholomorphic basis. In this case the polarized sections are simply the holomorphic functions. This example it is also called the holomorphic quantization.
\end{exmp}

\begin{exmp}[Real polarization]
Let $M=T^{\vee}C$ where $C$ is a configuration space. An example of real polarization is the vertical polarization associated to the cotangent bundle. The spann of $D$ is given by the momentum basis $\partial_{p_{j}}$. In this case the wave functions depends only by the position and the condition of preservation of $D$ on an observable $f$ is that:
	
$$ \frac{\partial^2}{\partial p_{j} \partial p_{k}}f\,=\,0 , \ \ \forall j,k=1, \ldots n.$$
\end{exmp}

The vertical polarization is tipically associated to the Sch\"{o}dinger representation of quantum mechanics. The holomorphic quantization brings to the Fock--Segal--Bergmann representation.

\subsection{K\"{a}hler quantization, holomorphic sections and Szeg\"{o} kernels}

In this subsection let us consider $(M,\omega,J,K)$ to be a compact K\"{a}hler manifold of complex dimension $n$. Assume we have the prequantization condition:

\begin{equation}
\label{PC5}
\Theta\,=\,-2i\omega.
\end{equation}

\begin{obs}
	We must do some considerations. The first is that the condition $(\ref{PC5})$ is slightly different from \textbf{PC4}, this is due by a different definition of $\omega$ but the geometrical essence is the same. Second, the condition $(\ref{PC5})$ ensures that the Dirac postulate is satisfied and, the fact that $\omega$ is integral, ensures the existence of the hermitian line bundle with $(\ref{PC5})$ satisfied.
\end{obs}

We call the triple $(L,\nabla,h)$ the prequantization bundle over $(M,\omega)$. We observe that with the hermitian product $h$ and the volume form $\omega^{n}=\omega\wedge \cdots \wedge \omega$, $n$ times defined on $M$, we can consider the space of smooth holomorphic sections $s$ of $M$ such that:

\begin{equation}
\label{vivahilbert}
\|s\|^2=\int_{M}h(s,s)\frac{\omega^{n}}{n!}
\end{equation}
\noindent
is finite. The $L^2$-completition of this space is an Hilbert space denoted by $\mathfrak{H}$ or also $H^0(M,L)$.

In the compact K\"{a}hler case from $(\ref{PC5})$ we have that $L$ is a positive line bundle and by the Kodaira embedding 
theorem there exists a positive tensor power $L^{\otimes k}$ with $k\in\mathbb{N}$ and global holomorphic sections $\{s^{k}_{i}\}_{i=0}^{n_{k}}$ that give the following embedding:

$$ \Phi: M \rightarrow \mathbb{P}^{n_{k}}(\mathbb{C}),$$
\noindent 
where $\Phi(z)=[s^{k}_{0}(z): \cdots :s^{k}_{n_{k}}(z)]$. The set $\{s^{k}_{i}\}_{i=0}^{n_{k}}$ is a basis for $H^{0}(M,L^{\otimes k})$ 
the space of holomorphic sections of $L^{\otimes k}$. The $\dimensione{(H^{0}(M,L^{\otimes k}))}=1+n_{k}$. 

In the contest of geometric quantization the parameter $k$ can be understood as a quantum parameter in the sense that $k=\frac{1}{\hbar}$. If we imagine that $\hbar\rightarrow 0$ then $k\rightarrow +\infty$ and we refind the semiclassical limit.
Let us consider the circle bundle $X$ of $L^{\vee}$ (defined previously) with $\pi:L^{\vee}\rightarrow M$ to $X$ that for simplicity we denote in the same way.
The circle bundle is the boundary of $D=\{(m,\lambda):m\in M,\lambda\in L_{m}^{\vee},h(\lambda,\lambda)\leq 1\} $ that is a strictly pseudoconvex domain in $L$.
We denote with $\|\cdot\|_{m}$ the induced norm by $h$ so we have that $D=\{\rho>0\}$ where $\rho:L^{\vee}\rightarrow \mathbb{R}$ 
is defined as $\rho(m,\lambda)=1-\|\lambda\|_{m}^{2}=1-a(m)^{-1}|\alpha|^2$ where we write $\lambda=\alpha s^{\vee}_{L}$ and $a(m)=\|s^{\vee}_{L}\|_{m}^{2}$ is a smooth function over a $U\subset M$ and 
$s^{\vee}_{L}$ is a local coframe over $U$. We have a circle action on $X$ denoted by 
$r_{\theta}:S^{1}\times X \rightarrow X$ with infinitesimal generator $\frac{\partial}{\partial \theta}$.   
As in \cite{Z} we consider the holomorphic and respectively antiholomorphic subspaces $T'D,T''D\subset TD_{\mathbb{C}}$ 
and the correspondent differentials $d'f=df|_{T'},d''f=df|_{T''}$ for $f$ smooth on $D$. $TD_{\mathbb{C}}=T'D\oplus T''D\oplus\mathbb{C}\frac{\partial}{\partial \theta}$
has a Cauchy Riemann structure and, the vectors on $D$ that are elements of $T'X$ (resp. $T''X$), are of the form $\sum_{j}a_{j}\frac{\partial}{\partial \overline{z}_{j}}$
(resp. $\sum_{j}a_{j}\frac{\partial}{\partial z_{j}}$). We can choose a basis for these vector spaces and consider the Cauchy Riemann operator
$\overline{\partial}_{b}:\mathcal{C}^{\infty}(X) \rightarrow C^{\infty}(X,(T''^{\vee}D))$ defined as $\overline{\partial}_{b}=df|_{T''}$.
If we define $\alpha=\frac{1}{i}d'\rho|_{X}$ and the Volume form $dV=\alpha\wedge(d\alpha)^{n}$ we have that $(X,\alpha)$ is a contact manifold.

\begin{defn}
	We define the Hardy space $H(X)=L^{2}(X)\cap \kernel(\overline{\partial}_{b})$ that admit the following decomposition:
\begin{equation}
\label{Hardy1}
H(X)=\bigoplus_{k}H(X)_{k} 
\end{equation}
\noindent 
where the subspaces 
\begin{equation}
\label{Hardy2}
H(X)_{k}=\{f\in C^{\infty}(X): f(e^{i\varphi}x)=e^{ik\varphi}f(x)\}\cap \kernel(\overline{\partial}_{b})
\end{equation}
\noindent	
are called $k$-Hardy spaces. 

\end{defn}

\begin{obs}
	On the Hardy spaces we have an hermitian product $H(s_{1},s_{2})=\int_{M}h_{m}(s_{1}(m),s_{2}(m))dV_{M}(m)=\int_{X}\widehat{s_{1}}\overline{\widehat{s_{2}}}dV_{X}$ because there is an unitary isometry 
	between sections $|\widehat{s}(x)|^2=\|s(m)\|^2$ (between $H^{0}(M,L^{\otimes k})$  and $H(X)_{k}$ with $m=\pi(x)$). The notation $\widehat{s}=\langle\lambda,s(z)\rangle$
    is used to denote the equivariant smooth section defined on $L^{\vee}\setminus\{0\}$. In analogue way $s_{k}$ on 
	$L^{\otimes k}$ determines a $\widehat{s}_{k}=\langle\lambda^{\otimes k},s_{k}(z)\rangle$. 
\end{obs}

\begin{defn}[Equivariant Szeg\"{o}  projector]
	We define the equivariant Szeg\"{o}  projector $\Pi_{k}:L^{2}(X) \rightarrow H(X)_{k}$ where $\forall f\in L^{2}(X)$ we have that: 
	
	$$\Pi_{k}(f)=\sum_{j}\langle f,\widehat{s}_{j}^{(k)}\rangle_{L^{2}(X)}\widehat{s}_{j}^{(k)}$$
\noindent
	where $(\widehat{s}_{j}^{(k)})_{j=1}^{n_{k}}$ is an orthonormal basis of $H(X)_{k}\cong H^{0}(M,L^{\otimes k})$.
\end{defn}

Expanding $\left(\Pi_{k}(f)\right)(x)=\int_{X}\sum_{j}\widehat{s}_{j}^{(k)}(x)\overline{\widehat{s}_{j}^{(k)}(y)}f(y)dV_{X}(y)$, we have this other definition.

\begin{defn}[Equivariant Szeg\"{o} kernel]
	The equivariant Szeg\"{o} kernel is:
	
	$$\Pi_{k}(x,y)=\sum_{j}\widehat{s}_{j}^{(k)}(x)\overline{\widehat{s}_{j}^{(k)}(y)}.$$
\end{defn}

By a theorem of \cite{BS} it is possible to represent the Szeg\"{o} kernel as a complex Fourier integral operator (FIO representation).

\begin{thm}
	Let $\Pi(x,y)$ be the Szeg\"{o} kernel of $X$, the boundary of a strictly pseudoconvex domain in $L$. Then there exists a symbol $s\in S^{n}(X\times X\times \mathbb{R}_{+})$ that admit the following 
	expansion:
	
	\begin{equation}
	\label{exp}
	s(x,y,t)=\sum_{k=0}^{\infty}t^{n-k}s_{k}(x,y)
	\end{equation}
	\noindent
	so that 
	
	\begin{equation}
	\label{fio}
	\Pi(x,y)=\int_{0}^{+\infty}e^{it\psi(x,y)}s(x,y,t)dt,
	\end{equation}
	\noindent 
	where $\psi\in C^{\infty}(D\times D)$ such that $\psi(x,x)=-i\rho(x)$ ($\rho$ define $D$), $d_{x}''\psi,d_{y}'\psi$ vanish to infinite order
	along the diagonal and $\psi(x,y)=-\overline{\psi(y,x)}$.
\end{thm}

We have that $\Pi$ is a Fourier integral operator with complex phase $\psi$ and the canonical relation $\Lambda \subset (T^{\vee}X) \times (T^{\vee}X)$ 
is generated by the phase $t\psi(x,y)$ on $X\times X\times\mathbb{R}_{+}$. In fact the canonical relation $\Lambda$ is the lagrangian submanifold of
$(T^{\vee}X) \times (T^{\vee}X)$ that has as generating function the phase function $t\psi$. 
The condition that must be true for the parametrization of the lagrangian submanifold $\Lambda$ is that:

\begin{equation}
\label{parametrizationC} \frac{d(t\psi)}{dt}=0
\end{equation}
\noindent 
that is when $\psi(x,y)=0$ and, on the diagonal $X\times X$ we have $ d_{x}{\psi}=-d_{y}\psi=-id'\rho|_{X}$. Let $\alpha=-id'{\rho}$ and let 
$ \Sigma=\{(x,r\alpha): r\in\mathbb{R}_{+}\}$ be the symplectic cone generated by the contact form $\alpha$, the real points of $\Lambda$ consist in the diagonal $\Sigma\times\Sigma$. We say that $\Pi$ has a Toeplitz structure on the symplectic cone $\Sigma$.

\begin{obs}
	The canonical relation $\Lambda$ cover an important rule in the theory of quantization as Guillemin and Sternberg remark in \cite{GS2} ``the smallest subsets of classical phase space in which the presence of a quantum mechanical particle can be detected are its lagrangian submanifolds''. 
\end{obs}

Now we are ready to recall an important result due to \cite{Z} that uses the method of stationary phase and the microlocal analysis of Szeg\"{o} and Bergman kernels.

\begin{thm}[Zelditch, Tian, Yau]
	Let $M$ a compact complex projective manifold of dimension $n$, and let $(A,h)$ a positive hermitian holomorphic line bundle. Let $g_{J}$ a K\"{a}hler metric on $M$ and 
	$-2i\omega=\Theta$ a K\"{a}hler form. For each $k\in\mathbb{N}$, $h$ induces a hermitian metric $h_{k}$ on $L^{\otimes k}$. Let $\{s^{k}_{i}\}_{i=0}^{n_{k}}$ be any 
	orthonormal basis of $H^{0}(M,L^{\otimes k})$ with $\dimensione{(H^{0}(M,L^{\otimes k}))}=1+n_{k}$. Then there exists a complete asymptotic expansion:
	
	\begin{equation}
	\label{ZZZ}
	\Pi_{k}(z,z)=\sum_{j}\|s_{j}^{(k)}(z)\|_{h_{k}}^2=a_{0}k^{n}+a_{1}(z)k^{n-1}+\cdots 
	\end{equation}
	
	for some $a_{j}$ smooth with $a_{0}=1$. 
\end{thm}

The proof is on \cite{Z} but we recall briefly the scheme. First we observe that $\Pi_{k}$ are Fourier coefficients of $\Pi$, so we have the following expression:

\begin{equation}
\Pi_{k}(x,y)=\int_{0}^{+\infty}\int_{S^{1}}e^{-ik\theta}e^{it\psi(r_{\theta}x,y,t)}s(r_{\theta}x,y,t)dtd\theta,
\end{equation}
\noindent 
where we used the FIO representation of $\Pi$ and $r_{\theta}$ denote the $S^{1}$ action on $X$. Changing variables $t\mapsto kt$ we have an oscillating integral:

\begin{equation}
\label{zelditchhhh}
\Pi_{k}(x,y)=k\int_{0}^{+\infty}\int_{S^{1}}e^{ik\Psi(t,\theta,x,y)}s(r_{\theta}x,y,t)dtd\theta
\end{equation}
\noindent
with phase $\Psi(t,\theta,x,y)=t\psi(r_{\theta}x,y)-\theta$. To simplify the phase we consider an holomorphic coframe $s_{L}^{\vee}$ and $a(z)=|s_{L}^{\vee}|_{z}^{-2}$. We write
any $\sigma\in L^{\vee}_{z}$ as $\sigma=\lambda s_{L}^{\vee}$ and for the coordinates $(x,y)=(z,\lambda,w,\mu)\in X\times X$ we have 
$\rho(z,\lambda)=a(z)|\lambda|^{2},\psi(z,\lambda,w,\mu)=-ia(z,w)\lambda\overline{\mu}$ with $a(z,w)$ is an almost analytic function on $M\times M$ such that $a(z,z)=a(z)$.
On $X$ we have $\lambda=a(z)^{-\frac{1}{2}}e^{i\phi}$  and $\mu=a(w)^{-\frac{1}{2}}e^{i\phi'}$. So for $(x,y)=(z,\phi,w,\phi')\in X\times X$ we have:

\begin{equation}
\label{phase}
\psi(z,\phi,w,\phi')=-i\left(1-\frac{a(z,w)}{i\sqrt{a(z)}\sqrt{a(w)}}\right)e^{i(\phi-\phi')}
\end{equation}
\noindent
and on the diagonal we have $ \psi(r_{\theta}x,x)=-i\left(1-\frac{a(z,z)}{a(z)}e^{i\theta}\right)$ and $\Psi(t,\theta,x)=-it(1-e^{i\theta})-\theta$.
The critical points for $\Psi$ are $\theta=0$ and $t=1$, the Hessian $\Psi''$ at this point is $\left(\begin{array}{cc} 0 & 1 \\ 1 & i\end{array}\right)$.
Applying the stationary phase method we find that:

$$ \Pi_{k}(x,x)\sim k\frac{\sqrt{2\pi i}}{\sqrt{\det{(k\Psi''(1,0))}}}\sum_{j,k=0}^{+\infty}k^{n-k-j}L_{j}s_{k}(x,x) $$
\noindent 
with $L_{j}$ differential operator of order $2j$ defined by:  

$$ L_{j}s_{k}(x,x)=\sum_{\nu-\mu=j}\sum_{2\nu\geq 3\mu}i^{-j}2^{-\nu}\langle\Psi''(1,0)^{-1}D,D\rangle^{\nu}\frac{(g_{(t,\theta)}^{\mu}ts_{k}(r_{\theta}x,x))|_{(t=1,\theta=0)}}{\mu!\nu!}$$
\noindent
where $\langle\Psi''(1,0)^{-1}D,D\rangle=2\frac{\partial^2}{\partial t\partial\theta}-i\frac{\partial^2}{\partial t^{2}}$.

The hypothesis of the theorem 7.7.5 of \cite{Ho} are satisfied  because the phase has nonnegative imaginary part and critical points are real and independent by $x$. 
So we have that $\Pi_{k}$ has the following form:

\begin{equation}
\Pi_{k}(x,x)=k^{n}C_{n}s_{0}(x,x)+k^{n-1}a_{1}(x,x)+\cdots
\end{equation}

The term $s_{0}(x,y)$ was expressed in \cite{BS} using that $\Pi$ is a projection and on the diagonal one has that
$s_{0}(x,x)dV=\frac{1}{4\pi^{n}}\det{L_{X}}\|d\rho\|dx$ with $L_{X}$ the restriction of the Levi form $L_{\rho}(z)=\sum\frac{\partial^2\rho}{\partial z_{j}\partial\overline{z}_{k}}z_{j}\overline{z}_{k}$ to the maximal complex subspace of $TX$. So:

\begin{equation}
\Pi_{k}(x,x)dV=k^{n}C_{n}\alpha\wedge\omega^{n}+O(k^{n-1}).
\end{equation}

The main coefficient is a nonzero constant $a_{0}$ times $k^{n}$. Comparing with the leading term of the Hirzebruch Riemann Roch polynomial we have that $a_{0}=1$. 
This conclude the proof. 

The method of the asymptotic expansion of the Szeg\"{o} kernel in presence of symmetries can be generalized in different ways. The interested reader can consult the following literature  \cite{P1}, \cite{P2},\cite{Cam}, \cite{GP1} and \cite{GP2}.

\section{The corrected geometric quantization: the case of complex cotangent bundle}

\subsection{The square--root bundle, the half form hilbert space and the harmonic oscillator}

From now we restrict our attention for the case where $M=T^{\vee}Q$ with $Q$ a configuration space of dimension $n$. Now $M=T^{\vee}Q$ is a K\"{a}hler manifold with K\"{a}hler form $\omega$, complex structure $J$ and K\"{a}hler potential $K$. We will consider associated to $M$ a K\"{a}hler polarization $D$.

In this subsection we recall some facts about determinant bundles, ``square--root bundles" and metaplectic correction.  

\begin{defn}
Let $D$ be a K\"{a}hler polarization of a K\"{a}hler manifold $(M,\omega,J,K)$. The deteminant bundle $\mathfrak{K}$ is the complex line bundle for which the sections are $n$--forms $\tau$ satisfying:

\begin{equation}
\label{determinant_b1}
\xi\,\lrcorner\, \tau \,=\, 0, \ \ \forall \xi \in\overline{D}.
\end{equation}

A section $\tau$ is said polarized if 

\begin{equation}
\label{determinant_b2}
\xi\,\lrcorner\, d\tau \,=\, 0, \ \ \forall \xi \in\overline{D}.
\end{equation}
\end{defn}

\begin{defn}
By square--root of $\mathfrak{K}$ we will denote the complex line bundle $\delta$ over $M$ such that $\delta \otimes \delta$ is isomorphic to $\mathfrak{K}$.
\end{defn}

On $\delta$ we have a partial connection $\nabla^{\delta}$ defined for $\xi\in \overline{D}$ and given by $\nabla_{\xi}\tau=\xi \,\lrcorner\, d\tau$. $\nabla^{\delta}$ descends from $\mathfrak{K}$ to $\delta$. So we have a connection of $L$ denoted by $\nabla^{L}$ and a partial connection on $\delta$ denoted by $\nabla^{\delta}$. In conclusion we have a partial connection on $L\otimes \delta$. 

\begin{prop}
If $\tau$ is an $(n,0)$--form  on $M$, then for each point $m$ in $M$ the $2n$--form:

\begin{equation}
\label{half_f1}
(-1)^{\frac{n(n-1)}{2}}(-i)^n\overline{\tau}\wedge \tau,
\end{equation}
\noindent 
is a non negative multiple of the volume form $dV_{M}\,=\, \frac{\omega^n}{n!}$. There is an unique hermitian structure on $\delta$ such that for each section $s$ of $\delta$:

\begin{equation}
\label{half_f2}
h(s,s)\,=\, \|s\|^2\,=\, \sqrt{\frac{(-1)^{\frac{n(n-1)}{2}}(-i)^n\overline{(s\otimes s)}\wedge (s\otimes s)}{2^n dV_{M}}}.
\end{equation}
\end{prop}

$\Proof.$\\
Let $z_{1}, \ldots , z_{n}$ be the holomorphic coordinates on $M$, then:
\begin{equation}
\label{pass1}
\omega^{n}\,=\, n!\det{(ig)}dz_{1}\wedge d\overline{z_{1}}\wedge \ldots \wedge dz_{n}\wedge d\overline{z_{n}},
\end{equation}
\noindent 
reordering differentials:

\begin{equation}
\label{pass1}
\frac{\omega^{n}}{n!}\,=\, \det{(ig)}(-1)^{\frac{n(n-1)}{2}}d\overline{z_{1}}\wedge\ldots \wedge d\overline{z_{n}}\wedge \ldots \wedge dz_{1}\wedge \ldots \wedge dz_{n},
\end{equation}
\noindent 
the second member is equivalent to $(\ref{half_f1})$ with $\overline{\tau}\,=\,d\overline{z_{1}}\wedge \ldots \wedge d\overline{z_{n}}$ and $\tau\,=\,dz_{1}\wedge\ldots \wedge dz_{n}$. 
\hfill $\Box$

Both the hermitian structure on $L$ and on $\delta$ gives in a natural way an hermitian structure on $L\otimes \delta$. 

\begin{defn}
	We define the half form Hilbert space $\mathfrak{H}^{L\otimes \delta}$ as the space of square integrable polarized sections of $L\otimes \delta$.
\end{defn}

Before to see how observables are quantizable in this new Hilbert space we recall some properties of sections of $L\otimes \delta$. First we can decompose a section $s$ of  $L\otimes \delta$ as $s\,=\, s^{L}\otimes s^{\delta}$. For $s^{\delta}$ we have that:

\begin{equation}
\label{derivative}
L_{\xi}(s^{\delta}\otimes s^{\delta})\,=\, s^{\delta}\otimes L_{\xi}s^{\delta}+L_{\xi}s^{\delta}\otimes s^{\delta} , \ \  \nabla_{\xi} (s^{\delta}\otimes s^{\delta})\,=\, 2\nabla_{\xi}(s^{\delta})\otimes s^{\delta},
\end{equation}
\noindent 
where $\xi$ is a vector field that preserve $\overline{D}$ and $s^{\delta}\otimes s^{\delta}\,=\, \left(s^{\delta}\right)^{2} \,=\, \tau\in\mathfrak{K}$. 

We have that polarized sections $s\,=\, s^{L}\otimes s^{\delta}$ of $L\otimes \delta$ are such that:

\begin{equation}
\label{polarized}
\nabla_{\xi}s \,=\, \nabla_{\xi}s^{L}\otimes s^{\delta} + s^{L}\otimes \nabla_{\xi}s^{\delta} \,=\, 0, \ \ \forall \xi \in \overline{D}.
\end{equation} 

Now if $f$ is an observable on $M$ such that $\xi_{f}$ preserves $\overline{D}$ then the quantum operator is given by:

\begin{equation}
\label{quantum_corrected1}
Q(f)\,=\, \left(Q_{\text{preq}}s^{L}\right)\otimes s^{\delta} - i\hbar s^{L}\otimes L_{X_{f}}(s^{\delta}),
\end{equation}
\noindent 
where $Q_{\text{preq}}$ is the prequantum operator and $s\,=\, s^{L}\otimes s^{\delta}$ is a section of $L\otimes \delta$. Note that:

\begin{equation}
\label{hamiltonian_v_f}
\xi_{f}\,=\, 2i(\partial_{\overline{z}}f\partial_{z}-\partial_{z}f\partial_{\overline{z}}),
\end{equation}
\noindent
then

$$ L_{\xi_{f}}s^{\delta}\,=\, L_{\xi_{f}}\sqrt{\tau} \,=\, \frac{1}{2}\tau^{-\frac{1}{2}}L_{\xi_{f}}\sqrt{\tau}\,=\, \frac{1}{2}\tau^{-\frac{1}{2}}L_{2i\partial_{\overline{z}}f\partial_{z}}\tau\,=\,  \frac{1}{2}\left(2i\partial_{z}\partial_{\overline{z}}f\right)\sqrt{\tau}. $$

If we look for $[\xi_{f},\partial_{\overline{z}}]\,=\,2i\partial_{\overline{z}}\partial_{z}f\partial_{\overline{z}}$ that is exactly the previous result. Then: 

\begin{equation}
\label{quantum_corrected1}
Q(f)(s^{L}\otimes\sqrt{\tau})\,=\, \left(Q_{\text{preq}}s^{L}\right)\otimes \sqrt{\tau} - i\hbar s^{L}\otimes\frac{1}{2}\left(2i\partial_{z}\partial_{\overline{z}}f\right)\sqrt{\tau}.
\end{equation}

\begin{exmp}
	In the case of $M=\mathbb{C}$, with hamiltonian $\mathcal{H}=\frac{1}{2}z\overline{z}$, the prequantum operator was:
	
	$$ Q_{\text{preq}}\,=\, \hbar z\partial_{z}$$
	\noindent 
	with eigenfunctions (sections) $\psi_{n}\,=\, z^n$ and eigenvalues $\hbar n$. These value of the energy are not the desired values. Let us consider a polarization $D$ spanned by $\partial_{\overline{z}}$. Then
	
	 $$2i\partial_{z}\partial_{\overline{z}}\mathcal{H}\,=\, i,$$
	 \noindent
	 and the corrected quantization gives:
	 
	 $$Q(f)\,=\, \hbar z \partial_{\overline{z}} +\frac{\hbar}{2}.$$
	 
	 This operators gives $\hbar \left(n + \frac{1}{2}\right)$, the correct eigenvalues. 
\end{exmp}

Details are in \cite{Ha} and \cite{W}.

\subsection{The BKS pairing}

Let us consider two Hilbert spaces $\mathfrak{H}_{1}$ and $\mathfrak{H}_{2}$ given by two different polarizations $D_{1}$ and $D_{2}$ on the same line bundle $L\rightarrow M$. The two Hilbert spaces are both part of the prequantum Hilbert space $\mathfrak{H}$. Let $P: \mathfrak{H}_{2}\rightarrow \mathfrak{H}_{1}$ be the orthogonal projection. Let $\langle \cdot ,\cdot \rangle$ be the inner product on $\mathfrak{H}_{1}$, we define the pairing $\langle\langle \cdot, \cdot \rangle \rangle: \mathfrak{H}_{1}\times \mathfrak{H}_{2}\rightarrow \mathbb{C}$ by:

\begin{equation}
\label{pairing}
\langle\langle s_{1},s_{2}  \rangle \rangle\,=\, \langle s_{1},Ps_{2}\rangle.
\end{equation}

The map is not unitary but if the flow of an observable $f\in C^{\infty}(M)$ preserves both polarizations then:

\begin{equation}
\label{pairing2}
\langle\langle Q(f)s_{1},s_{2}  \rangle \rangle\,=\, \langle\langle s_{1},Q(f)s_{2}  \rangle \rangle.
\end{equation}

Note that when $D_{1}$ and $D_{2}$ are determinated by two complex structures on the symplectic space then a multiple of the projector is in fact unitary. 

\subsection{The BKS pairing on the cotangent bundle, the time evolution and the Schr\"{o}dinger equation}

In this subsection we enter in details examining the case of cotangent bundles $M=T^{\vee}C$, for some configuration space $C$. In this case $M$ is a model for classical mechanics, a canonical transformation $\rho$ preserves the symplectic structure $\omega$. We can consider an Hamiltonian vector field $X_{f}$ that generates a canonical flow $\rho_{t}$. This flow induces a time evolution of observable via $f_{t}(m)=f(\rho_{t}(m))$ (pull-back for every $m\in M$). Now this time evolution is lifted, through $Q_{\text{preq}}(f)$, at the level of sections of an hermitian line bundle $L\rightarrow M$. In particular the lifting flow is: 

\begin{equation}
\label{lift_flow_quantum}
\hat{\rho}_{t}s(m)\,=\, e^{-\frac{1}{\hbar}\int_{0}^{t}\mathcal{L}ds}s(\rho_{t}(m)),
\end{equation}
\noindent 
where $\mathcal{L}\,=\, X_{f}\,\lrcorner\,\theta -f$ (details  are in \cite{C}) and this lift is unitary. Given an observable $f\in C^{\infty}(M)$ that preserve a polarization $D$ then we describe the unitary time evolution by:

\begin{equation}
\label{lift_flow_quantum2}
Q_{\text{preq}}(f)s_{t}\,=\, \frac{d}{dt'}(\hat{\rho}_{t'-t}s_{t})|_{t'=t}\,=\, -i\hbar\dot{s}_{t}.
\end{equation}

Considering the half--form construction, that is the line bundle $L\otimes \delta$, we use a pull-back in order to define an evolution in $t$:

\begin{equation}
\label{lift_flow_quantum2}
\rho_{t}^{L\otimes \delta}s_{t}^{L\otimes \delta}\,=\, (\hat{\rho}_{t}s^{L})(s^{\delta}(\rho_{t})).
\end{equation}

The derivative at $t=0$ reproduce the quantum operator corrected. Now let us assume on $M$ two transverse real polarizations $D_{1}$ and $D_{2}$. Then $TM=D_{1}\oplus D_{2}$. We can write $M$ as the product of two configuration spaces $C\times C'$, one for position $q$ and the second as the momentum space $q'$. Then we have that:

\begin{equation}
\label{canonical_omega}\omega=\partial_{q}\partial_{q'}\mathcal{S}dq\wedge dq',
\end{equation}
\noindent 
with $\mathcal{S}$ the generating function. We can define a pairing between $\delta_{D_{1}}$ and $\delta_{D_{2}}$:

\begin{equation}
\label{pairing_half_forms}
(\sqrt{\tau_{1}},\sqrt{\tau_{2}})\,=\, \sqrt{\tau_{2}\wedge \overline{\tau_{1}}},
\end{equation}
\noindent 
where $\tau_{1},\tau_{2}$ are $n$-form of the determinant bundle. Then the BKS pairing between sections of $\mathfrak{H}_{1}$ and $\mathfrak{H}_{2}$ is:

\begin{equation}
\label{particular_pairing}
\langle\langle s_{1}^{L\otimes\delta_{D_{1}}},s_{2}^{L\otimes\delta_{D_{2}}}\rangle\rangle\,=\, \int_{M}(s_{1}^{L},s_{2}^{L})(\sqrt{\tau_{1}},\sqrt{\tau_{2}})dV_{M}.
\end{equation}

The Liouville measure is $dV_{M}=\frac{\omega^{n}}{(2\pi \hbar)^n}=\frac{\det{\partial_{q}\partial_{q'}\mathcal{S}}}{(2\pi \hbar)^n}(dq\wedge dq ')^n$ and a computation shows that $(\sqrt{\tau_{1}},\sqrt{\tau_{2}})dV_{M}\,=\, \frac{\sqrt{\det{\partial_{q}\partial_{q'}\mathcal{S}}}}{(2\pi \hbar)^{n/2}}d^nq\wedge d^nq'$.

Let $\mathcal{S}=p\cdot q$ be the canonical generator, then the determinant will be the identity, $q'=p$ corresponds to the base space of $D_{2}$ whose leaves are surfaces of constant $q'=p$ and, because the polarization is transverse to $D_{1}$ the surfaces of constant $q$ are leaves of $D_{1}$. At the end the pairing is given by:

\begin{equation}
\label{fourier1}
\langle\langle s_{1}^{L\otimes\delta_{D_{1}}},s_{2}^{L\otimes\delta_{D_{2}}}\rangle\rangle\,=\,\frac{1}{(2\pi\hbar)^{n/2}}\int_{C}\overline{\psi}(q)\int_{C'}\psi(p)e^{\frac{ip\cdot q}{\hbar}}d^npd^nq,
\end{equation}
\noindent 
where $\phi(q)$ and $\psi(p)e^{\frac{i}{\hbar}p\cdot q}$ are complex wave functions. The projector corresponds effectively to the Fourier transform:

\begin{equation}
\label{fourier2}
Ps_{2}\,=\,\frac{1}{(2\pi\hbar)^{n/2}}\int_{C'}\psi(p)e^{\frac{ip\cdot q}{\hbar}}d^np,
\end{equation}
\noindent 
with value in $\mathfrak{H}_{D_{1}}$.

\subsection{Reconstruction of Schr\"{o}dinger equation on the cotangent bundle}

Let $s^{L\otimes \delta}=s^{L}\otimes \sqrt{\tau}$ be a $D_{1}$--polarized section, $X_{f}$ an hamiltonian vector field and $\hat{\rho}_{t}s^{L\otimes \delta}=\hat{\rho}_{t}s^{L}\rho_{t}^{*}\sqrt{\tau}$ the lift of the canonical flow in $D_{2}=\rho^{*}_{t}D_{1}$. Then:

\begin{equation}
\label{dynamics1}
\langle \dot{s}^{L\otimes \delta},w^{L\otimes \delta}\rangle \,=\, -\frac{d}{dt'}\left. \langle \langle \hat{\rho}_{t'-t}s^{L\otimes \delta}_{t},w\rangle\rangle \right|_{t'=0},
\end{equation}
\noindent 
this is for every $w\in\mathfrak{H}(D_{1})$ where $w=w^{L}\sqrt{\sigma}$ and for $s^{L\otimes \delta}_{t}$ we assume equal to $s^{L\otimes \delta}_{0}$. We have that:

\begin{equation}
\label{dynamics2}
\langle \dot{s}^{L\otimes \delta},w^{L\otimes \delta}\rangle \,=\, -\lim_{h\rightarrow 0}\frac{\langle P\hat{\rho}_{h-t}s^{L\otimes \delta}_{t},w\rangle-\langle P\hat{\rho}_{-t}s^{L\otimes \delta}_{t},w\rangle}{h}.
\end{equation}

We must to evaluate the following expression:

\begin{equation}
\label{dynamics3}
\langle \langle \hat{\rho}_{t}s^{L\otimes \delta},w^{L\otimes \delta}\rangle\rangle \,=\, \int_{M}\overline{\psi(\rho_{t}(m))}w(q)e^{\frac{i}{\hbar}\int_{0}^{t}\mathcal{L}(\gamma(s))ds}(\rho_{t}^{*}\sqrt{\sigma},\sqrt{\sigma})dV_{M}.
\end{equation}

In our situation we consider the free particle on $M=\mathbb{R}^{2n}$, with $\mathcal{L}=\mathcal{H}=\frac{p^2}{2m}$, $\rho_{t}$ the uniform motion on flat space (curvature term is $0$) with equation: $q(t)=q+vt=q+\frac{p}{m}t$ and $p(t)=p$. We observe that: $\int_{0}^{t}\mathcal{L}(\gamma(s))ds=\mathcal{H}t=\frac{p^2t}{2m}$. Let $x:M\rightarrow \mathbb{R}^{2n}$ be the coordinates with $x(\rho_{t}(m))=(q(t),p(t))$ and $dq^{n}=dq_{1}\wedge \ldots \wedge dq_{n}$. Then:

$$ (\rho_{t}^{*}\sqrt{\sigma},\sqrt{\sigma})^2dV_{M} \,=\, d(x\circ \rho_{t}) \,=\, d\left(q+\frac{t}{m}p\right)^n\wedge dq^n\left(\frac{1}{(2\pi\hbar)^n}\right)=\left(\frac{t}{2\pi m \hbar}\right)^n dp ^{n}\wedge dq^n.$$

That is: 

$$ (\rho_{t}^{*}\sqrt{\sigma},\sqrt{\sigma})dV_{M}\,=\,\left(\frac{t}{2\pi m \hbar}\right)^{\frac{n}{2}} dp ^{n}\wedge dq^n.$$

After these refinements we must to evaluate the integrals:

\begin{equation}
\label{dynamics4}
\langle \langle \hat{\rho}_{t}s^{L\otimes \delta},w^{L\otimes \delta}\rangle\rangle \,=\, \int_{\mathbb{R}^{2n}}\overline{\psi(q(t),p(t))}w(q)e^{\frac{i}{\hbar}\frac{p^2t}{2m}}\left(\frac{t}{2\pi m \hbar}\right)^{\frac{n}{2}} dp^{n}\wedge dq^n.
\end{equation}

In order to evaluate the integral we start as in \cite{C} to expand the function $\psi$ around the point $q$.

\begin{equation}
\label{taylor_exp_around_q}
\begin{multlined}[12cm]
\psi(q(t),p(t))\,=\, \psi\left(q+\frac{p}{m}t,p\right)\,=\, \psi(q)+ \\ +\frac{t}{m}\sum_{j}p_{j}\partial_{q_{j}}\psi(q)+ \frac{t^2}{2m^2}\sum_{j,l}p_{j}p_{l}\partial_{q_{j}}\partial_{q_{l}}\psi(q)+ \mathcal{O}(t^3).
\end{multlined}
\end{equation}

The contribution given by the second term in the expansion to the initial integral is $0$ due by the parity of the integrand. Substituting inside the integral:

\begin{equation}
 \label{dynamics5}
 \int_{\mathbb{R}^{n}}w(q)\left(\frac{t}{2\pi m \hbar}\right)^{\frac{n}{2}}\int_{\mathbb{R}^{n}}\left[\overline{\psi}(q)+ \frac{t^2}{2m^2}\sum_{j,l}p_{j}p_{l}\partial_{q_{j}}\partial_{q_{l}}\overline{\psi}(q)\right]e^{\frac{i}{\hbar}\frac{p^2t}{2m}} dp^{n} dq^n + \mathcal{O}(t^3).
 \end{equation}

Now we can evaluate the two integral respect the momentum variable $p$. The first is:

\begin{equation}
\label{first_integral_Schrodinger}
\int_{\mathbb{R}^{n}}e^{\frac{i}{2}\frac{t}{m\hbar}p^2}dp^n \,=\, \frac{(2\pi)^{\frac{n}{2}}}{\left(\frac{t}{m\hbar}\right)^{\frac{n}{2}}}e^{ic\frac{\pi}{4}},
\end{equation}
\noindent 
where $e^{ic\frac{\pi}{4}}$ is a phase that depends by the metaplectic correction. This integral as been estimated using the stationary phase method on multidimensional Fresnel integrals (the reference is \cite{S}). 
The second integral to estimate is the following:

\begin{equation}
\label{second_integral_Schrodinger}
\int_{\mathbb{R}^{n}}p_{j}p_{l}e^{\frac{i}{2}\frac{t}{m\hbar}p^2}dp^n \,=\,\left(\int_{\mathbb{R}}p_{j}^2e^{\frac{i}{2}\frac{t}{m\hbar}p_{j}^2}dp_{j}\right)\cdot\left(\int_{\mathbb{R}^{n-1}}e^{\frac{i}{2}\frac{t}{m\hbar}\widehat{p}^2}d\widehat{p}^{n-1}\right),
\end{equation}
\noindent 
where the only terms different from zero inside the integral are when $j=l$, so we can consider this integral as one dimensional integral with amplitude $p_{j}^{2}$ (evaluating using the result of \cite{AM}) and a second $n-1$ integral as before where $\widehat{p}$ is obtained by $p$ excluding the $j$ coordinate. We find that:

\begin{equation}
\label{second_integral_Schrodinger2}
\int_{\mathbb{R}^{n}}p_{j}p_{l}e^{\frac{i}{2}\frac{t}{m\hbar}p^2}dp^n \,=\,e^{i\frac{\pi}{4}}\left(\frac{2m\hbar}{t}\right)^{\frac{1}{2}}\left[\frac{e^{i\frac{\pi}{2}}}{2}\left(\frac{2m\hbar}{t}\right)\frac{\sqrt{\pi}}{2}\cdot 2\right]\cdot \frac{(2\pi)^{\frac{n-1}{2}}}{\left(\frac{t}{m\hbar}\right)^{\frac{n-1}{2}}}e^{i\widetilde{c}\frac{\pi}{4}},
\end{equation}
\noindent 
where $c=\widetilde{c}+1$. After calculations and semplifications the integral $(\ref{dynamics5})$ assume the following form:

\begin{equation}
\label{dynamics6}
\int_{\mathbb{R}^{n}}w(q)e^{ic\frac{\pi}{4}}\left[\overline{\psi}(q)+ \frac{it\hbar}{2m}\nabla\overline{\psi}(q)\right]dq^n + \mathcal{O}(t^2),
\end{equation}
\noindent 
that we can write as:

\begin{equation}
\label{dynamics7}
\int_{\mathbb{R}^{n}}w(q)e^{ic\frac{\pi}{4}}\left[\overline{\psi}(q)- \frac{it}{\hbar}\mathcal{H}\overline{\psi}(q)\right]dq^n + \mathcal{O}(t^2),
\end{equation}
\noindent 
with $\mathcal{H}\,=\, -\frac{\hbar^2}{2m}\nabla^2$. Now after a differentiation respect to $t$ at $t=0$ and taking the complex conjugate, we find for every section of $L\otimes \delta$:

\begin{equation}
\label{dynamics8}
i\hbar\frac{\partial}{\partial t}\psi(q)\,=\,\mathcal{H}\psi(q) ,
\end{equation}
\noindent 
the Schr\"{o}dinger equation. 

\subsection{The Klein--Gordon equation and the cotangent bundle}

Let $Y$ be a space--time manifold. Assuming that $Y$ is an orientable manifold, we have that $T^{\vee}Y$ is a symplectic manifold with symplectic form:

\begin{equation}
\label{emf_sympl}\omega_{e}\,=\, d\theta_{Y} + e\pi^*F,
\end{equation}
\noindent 
where $F$ is the electromagnetic field (a closed 2-form on $Y$), $\theta_{Y}$ is the symplectic potential,  $e$ is the charge of a particle in an external electromagnetic field and $\pi:T^{\vee}Y\rightarrow Y$.  

The metric $g$ is the Lorentz metric with signature $(-,-,-,+)$. Let us consider the function $N(x)\,=\, g(x,x)$ defined on $T^{\vee}Y$. The square root of $N$ represents the mass of a particle in the classical state $x$. Let us consider the complete hamiltonian vector field $\xi_{N}$ and the flux $\phi_{t}^{N}$ generated by $\xi_{N}$. 

In this case the prequantum condition reads as:

$$ \int_{\Sigma}eF \in h\mathbb{Z},$$
\noindent 
for every compact oriented 2-surface $\Sigma$ on $Y$ with empty boundary. 

We can consider a vertical polarization $D^{\mathbb{C}}$ on $T^{\vee}Y$ and the associated Hilbert space $\mathfrak{H}_{D^{\mathbb{C}}}$ that is the $L^{2}(Y,g)$. 

Under the assumption that $D^{\mathbb{C}}$  and $T\phi_{t}^{N}(D^{\mathbb{C}})$ are transverse, it is possible to evaluate $Q_{\text{preq}}(N)s_{t}$. The method is similar as before and is showed in \cite{Sn} in the chapter 10. Let us denote with $L_{e}$ the corrected quantum line bundle, then the BKS kernel $BKS_{t}$ is given by:

\begin{equation}
\label{BKS_relativistic}
BKS_{t}( s_{1},s_{2} )\,=\, \int_{T^{\vee}Y}\langle s_{1},\phi_{t}^{N}s_{2}\rangle \left|\omega_{e}^{4}\right|,
\end{equation}
\noindent 
for local sections $s_{1}\,=\, \psi_{1}s_{1}^{L}\otimes s_{1}^{\delta} , s_{2}\,=\, \psi_{2}s_{2}^{L}\otimes s_{2}^{\delta}$. 

What the author find at the end of calculations (in \cite{Sn}), is the following expression for the quantum operator:

\begin{equation}
\label{lift_flow_quantum2_final}
Q_{\text{preq}}(N)s\,=\, -\hbar^2\left[\Box-\frac{1}{6}R\right]s,
\end{equation}
\noindent 
where $\Box$ is the d'Alembert operator and $R$ is the scalar curvature of the metric $g$. The operator $N$ is also called the ``mass-squared operator'' and, the equation: 

$$-\Box\psi+\frac{1}{6}R\psi\,=\,0,$$
\noindent 
is the Klein-Gordon equation and $\psi$ the wave function. 

In \cite{BO} has been discovered that the critical metrics of the expectation value of $Q(N)$ have to satisfy Einstein's equation for suitable energy--momentum tensor.

\subsection{Considerations on the BKS pairing}

If the first problem of the prequantization procedure was to have the dependence on half variables by the wave functions (sections of the line bundle), solved by introducing the concept of polarizations, the second problem was the treatment of second order, or higher operators, as the free particle and harmonic oscillator. This second problem has been solved using the BKS construction with the half--forms bundle and the lift of the flow $\rho_{t}$ that describes the dynamic of the system. As focused also by \cite{C} this seems the ``exact" choice of the quantum operators, with a suitable metaplectic correction and considering the Hamiltonian dynamics  (so the presence of symmetries). The final form of our observables $f\in C^{\infty}(M)$, in the particular case of cotangent bundles for section $s\sqrt{\tau}$ in the corrected bundle, seems to be:

\begin{equation}
\label{lift_flow_quantum2_final}
Q_{\text{preq}}(f)s_{t}\,=\, -i\hbar\frac{\partial}{\partial t}P\hat{\rho}_{t}s|_{t=0},
\end{equation}
\noindent 
where $P$ is the usual Fourier transform. As observed always in \cite{C} generalization on operators of order bigger then $2$ is a problem. Generally the Dirac axiom is not true, the time evolution can be not unitary, operators are not always self--adjoint and operators can be not linear. These problems have been showed not only in  \cite{C} but also in the deep study of \cite{AE}. 

\subsection{The BKS pairing for one complex and one real polarization: the Segal--Bargmann transform}

Let us consider for $M$ a K\"{a}hler manifold and two polarizations: one $D_{1}$ a real polarization and $D_{2}$ a K\"{a}hler polarization. Let us assume $M=\mathbb{C}$, as in the example of the harmonic oscillator. We have that $D_{2}$ is spanned by $\partial_{\overline{z}}$ and $D_{1}$ corresponds to $q=$ constant. Let us denote the two Hilbert spaces by $\mathfrak{H}(D_{1})$ and $\mathfrak{H}(D_{2})$ and by $z=p+iq$ and $\overline{z}=p-iq$ the respective coordinates. The following symplectic potential:

\begin{equation}
\label{symplectic_CRpolarization_ex}
\theta\,=\, \frac{i}{2}z d\overline{z} \,=\, \frac{pdq-qdp}{2},
\end{equation} 
\noindent 
is associated to the K\"{a}hler form: $\omega=d\theta =dp\wedge dq$. The line bundle $L$ over $M$ is trivial and the space $\mathfrak{H}(D_{2})$ corresponds to the Fock space $\mathfrak{F}$ of holomorphic function with scalar product:

\begin{equation}
\label{CRpolarization_fock_space_product}
(f,g)\,=\, \frac{1}{2\pi}\int_{\mathbb{C}}f(z)\overline{g(z)}e^{-\frac{|z|^2}{2}}dzd\overline{z}.
\end{equation}

A section $s$ of the line bundle has squared norm $\|s\|^2\,=\, e^{-\frac{|z|^2}{2}}$ and $s$ have the following form:

\begin{equation}
\label{section_d2_CRpolarization_ex}
\phi(z)e^{-\frac{z\cdot \overline{z}}{4}},
\end{equation} 
\noindent 
where $\phi(z)$ is an holomorphic function. A section for the real polarization $D_{1}$ is the complex function:

\begin{equation}
\label{section_d1_CRpolarization_ex}
\psi(z)e^{-i\frac{pq}{2}}.
\end{equation} 

Let us consider the half form bundle $\delta$ in order to ensure the metaplectic correction. We can define $P:\mathfrak{H}(D_{2}) \rightarrow \mathfrak{H}(D_{1})$ and $P':\mathfrak{H}(D_{1})\rightarrow \mathfrak{H}(D_{2})$ the projections given by the BKS pairing. 
The pairing between the two hilbert spaces is given by the following formula:

\begin{equation}
\label{pairing__CRpolarization_ex}
\langle\langle \psi,\phi \rangle \rangle\,=\, \frac{C}{2\pi}\int_{\mathbb{R}^2}\overline{\psi}(q)\phi(z)e^{\frac{2ipq - p^2 - q^2}{4}}dpdq,
\end{equation}
\noindent 
where $C$ is a constant we specify at the end. Now we can give and explicit form for $P$ and $P'$ using the reproducing property of holomorphic function in the second hilbert space:

\begin{equation}
\label{reproducing_CRpolarization_ex}
\phi(z)\,=\, \frac{1}{2\pi}\int_{\mathbb{C}}\phi(w)e^{\frac{\overline{w}z-w\overline{w}}{2}}dwd \overline{w},
\end{equation}
\noindent 
where $w=x+iy$.

Now inserting this expression inside the BKS pairing $(\ref{pairing__CRpolarization_ex})$, we can integrate respect $p$ on $\mathbb{R}$ the expression in order to have:

\begin{equation}
\label{pairing__CRpolarization_ex2}
\langle\langle \psi,\phi \rangle \rangle\,=\, \frac{1}{2\pi}\int_{\mathbb{C}\times\mathbb{R}}\overline{\psi}(q)\phi(w)\overline{K}(q,w)e^{\frac{- w\overline{w}}{4}}dqdwd\overline{w},
\end{equation}
\noindent 
with

\begin{equation}
\label{pairing__CRpolarization_ex2_K}
K(q,w)\,=\, \frac{\overline{C}}{\sqrt{\pi}}e^{-\frac{q^2}{2}-iqw+\frac{w^2}{4}}.
\end{equation}

Then the expression for the two projections is:

\begin{equation}
\label{CRpolarization_ex2_P1}
P(\phi)(q)\,=\, \frac{1}{2\pi}\int_{\mathbb{C}}\phi(w)\overline{K}(q,w)e^{\frac{- w\overline{w}}{2}}dwd\overline{w},
\end{equation}
\noindent 
and

\begin{equation}
\label{CRpolarization_ex2_P1}
P'(\psi)(w)\,=\, \frac{1}{2\pi}\int_{\mathbb{R}}\psi(q)\overline{K}(q,w)dq.
\end{equation}

We observe that functions $z^m$ in $\mathfrak{H}(D_{2})$ are mapped in eigenfunctions of the harmonic oscillator. This is the Segal--Bargmann transform and the ground state $1$ is mapped to:

\begin{equation}
\label{CRpolarization_ex2_P1_armonic_oscillator}
P(\phi)(q)\,=\,\frac{C}{\sqrt{\pi}}e^{-\frac{q^2}{2}},
\end{equation}
\noindent 
with $C=\pi^{\frac{1}{4}}$.

\subsection{The BKS pairing between two complex polarizations: the Bogoliubov transformation}

The last case is the case of two complex polarizations $D_{1}$ and $D_{2}$ determined respectively by two complex structures $J_{1}, J_{2}$. In this case a multiple of the projector operator is unitary and the existence is guarantee by the Stone Von Neumann theorem (see \cite{W}). In this brief summary we consider the case in  one dimension with $M=\mathbb{C}$ and without metaplectic correction. To the polarization $D_{1}$ we have the Hilbert space $\mathfrak{H}(D_{1})$ that is the Fock space with holomorphic sections:

\begin{equation}
\label{Bogoliubov_sections1}
\phi(z)\,=\, \psi(z)e^{-\frac{1}{4}z\overline{z}},
\end{equation}
\noindent 
where $\psi(z)$ is holomorphic respect $J_{1}$ and $z=p+iq$. The inner product is:

\begin{equation}
\label{Bogoliubov_sections1}
(\phi(z),\phi'(z))\,=\, \frac{1}{2\pi}\int_{\mathbb{C}}\psi(z)\overline{\psi'}(z)e^{-\frac{1}{2}z\overline{z}}dzd\overline{z}.
\end{equation}

Now let $J_{2}$ be a second complex structure that determines $D_{2}$ and a second hilbert space $\mathfrak{H}(D_{2})$ (always a Fock space) with sections:

\begin{equation}
\label{Bogoliubov_sections2}
\phi'(w)\,=\, \psi'(w)e^{-\frac{1}{4}w\overline{w}},
\end{equation}
\noindent 
where $\psi'(z)$ is holomorphic respect $J_{2}$ (the inner product is the same of the previous case). The BKS pairing is given by:

\begin{equation}
\label{Bogoliubov_pairing3}
\langle\langle \phi(z),\phi'(z) \rangle\rangle\,=\, \frac{1}{2\pi}\int_{\mathbb{C}^2}\psi(z)\overline{\psi'}(w)e^{-\frac{1}{4}\left(z\overline{z}+w\overline{w}\right)}dzd\overline{z}dwd\overline{w}.
\end{equation}

The pairing determines a projection $P:\mathfrak{H}(D_{2})\rightarrow \mathfrak{H}(D_{1})$ given by:

\begin{equation}
\label{Bogoliubov_pairing4}
P\psi'(w)\,=\, \frac{1}{2\pi}\int_{\mathbb{C}}\overline{\psi}_{z}\psi'(w)dwd\overline{w},
\end{equation}
\noindent 
where $\psi_{z}(w)\,=\, e^{\frac{1}{2}(2\omega(z,(J_{1}+i)w)-\omega(z,J_{1}z)-\omega(w,J_{1}w))}$. The projection of $\phi_{0}'(w)=1\cdot e^{\frac{1}{4}w\overline{w}}$ (the ground state in $\mathfrak{H}(D_{1})$) is:

\begin{equation}
\label{Bogoliubov_pairing5}
\phi_{0}(z)\,=\,\left[\det{\frac{1}{2}(J_{1}+J_{2})}\right]^{-\frac{1}{2}}e^{\frac{\lambda(z)}{4}-\frac{z\overline{z}}{4}},
\end{equation}
\noindent 
where $\lambda(z)=2\omega(z,J_{1}Lz)-2i\omega(z,Lz)$ and $L=(J_{1}+J_{2})^{-1}\cdot(J_{1}-J_{2})$. A good reference is \cite{W2} where also the infinite dimensional case is explained. The same topic is present in \cite{W}. 

\section{The Bohr--Sommerfeld (sub)manifolds}

Let us consider a lagrangian submanifold $\Lambda$ of a cotangent space $T^{\vee}Q$. It is possible to interpret the quantity $e^{i\frac{z}{\hbar}}$ as a covariantly constant section $s$ of the line bundle $L_{\Lambda}$. Now there is a natural definition of square root of a volume form $\tau$ on $\Lambda$ replacing $\tau$ by $\sqrt{\tau}$. The section $\sqrt{\tau}$ is a section of $\sqrt{\Delta_{\Lambda}}$ where $\Delta_{\Lambda}$ is the line bundle of the complex volume forms on $\Lambda$. We can consider the triple $\left(\Lambda,s,\sqrt{\tau}\right)$ and to work with the machinery of geometric quantization. Following \cite{I}, we can study the asymptotic expansion of the squared norm of an ``isotropic state'' in particular conditions. An isotropic state is a family of sections of an Hilbert space associated to the lagrangian submanifold. Another problem is the study of finding  a simultaneous WKB eigenstate of a set of classical observables $\mathcal{H}_{1}, \ldots , \mathcal{H}_{n}$ in involution.

Assuming that $\Lambda$ must be of the form $h_{i}=E_{i}$ with $E_{i}$ that are constant for every $i=1, \ldots , n$. We can consider the volume form $\tau$ of $\Lambda$. Assuming for simplicity that $\Lambda$ is simply connected, the wave functions are well defined if $s\sqrt{\tau}$ is a single valued global section of $L_{\Lambda}\otimes \sqrt{\Delta_{\Lambda}}$. The existences of this section depends by the validity of the Bohr-Sommerfeld condition \textbf{PC5}. Assume that $L=M\times \mathbb{C}$ and that $\theta$ is the symplectic potential, then the corrected quantization condition is that:

\begin{equation}
\label{bcondition}
\left(e^{\frac{2i}{\hbar}\oint_{\gamma} \theta}\right)\tau
\end{equation}  
\noindent
a section of $\Delta_{L}$, admits a single values square root with values in $\sqrt{\Delta_{L}}$. Here $\gamma$ is any loops in $\Lambda_{m}$ ($m\in M$). There are other formulations of \textbf{PC5}, for example in \cite{B} it is given by:

\begin{equation}
\label{bcondition_2}
\oint_{\gamma}\theta\,=\, 2\pi\hbar(n+d),
\end{equation}  
\noindent
where $n\in \mathbb{Z}$, as in the first quantization condition \textbf{PC1}, with the difference that there is the quantity $d$ associated to the holonomy from the flat connection of the bundle of half-forms. 
Authors who worked in the setting of geometric quantization of Bohr--Sommerfeld (sub)manifolds are \cite{BPU},\cite{D} and \cite{DP} (where the last two concentrated the attention on the equivariant case).

\section{Conclusion}
The real contribution of this paper consists of a new way to derive the Sch\"{o}dinger equation, as in \cite{C}, with the difference we used a result of Albeverio and Mazzucchi (\cite{AM}) on asymptotic expansion of Fresnel integrals instead the heat kernel. All possible expressions of quantization conditions have been examined in detail. These notes provide a complete, detailed summary on the state of the art on the geometric quantization showing possible connection between mathematics and physics.

\end{document}